\documentclass{amsart}

\usepackage[utf8x]{inputenc}

\usepackage{ae,aecompl}

\usepackage{subfigure}

\usepackage{amssymb}
\usepackage{amsmath}
\usepackage{amscd}
\usepackage{units}

\usepackage[english]{babel}

\usepackage[normalem]{ulem}

\usepackage{tikz,ifthen}
\usetikzlibrary{matrix,shapes,arrows,chains}
\usetikzlibrary{calc}
\usetikzlibrary{patterns}
\usetikzlibrary{positioning} 

\definecolor{darkgreen}{rgb}{0,0.7,0}
\definecolor{darkred}{rgb}{0.7,0,0}
\definecolor{darkblue}{rgb}{0,0,0.7}
\usepackage{pdflscape}

\usepackage{enumerate}
\usepackage{verbatim}
\usepackage{hyperref}

\usepackage{graphicx}
\usepackage{rotating}
\usepackage{algorithmic}
\usepackage{algorithm}
\usepackage{xspace}

\usepackage{ifpdf}
\ifpdf
\fi

\newtheorem{theorem}{Theorem}[section]

\newtheorem{proposition}[theorem]{Proposition}
\newtheorem{corollary}[theorem]{Corollary}
\newtheorem{definition}[theorem]{Definition}

\newtheorem{conjecture}[theorem]{Conjecture}

\newtheorem{remark}[theorem]{Remark}

\def\SG{{\mathfrak S}}  
\def\x{{\mathbf x}}

\definecolor{darkgreen}{rgb}{0,0.6,0}
\definecolor{darkred}{rgb}{0.7,0,0}
\definecolor{darkblue}{rgb}{0,0,0.7}

\newboolean{draft}
\setboolean{draft}{true}
\ifdraft
\newcommand{\TODO}[2][To do: ]{\textcolor{red}{\textbf{#1#2}}}
\newcommand{\INFO}[2][Info: ]{\textcolor{darkgreen}{\textbf{#1#2}}}
\else
\newcommand{\TODO}[2][]{}
\newcommand{\INFO}[2][]{}
\fi

\author{Nicolas Borie}

\title{The Hopf Algebra of graph invariants}

\address{Univ. Paris Est Marne-La-Vall\'ee, Laboratoire d'Informatique
  Gaspard Monge, 5, bd Descartes Champs sur Marne 77454
  Marne-la-Vall\'ee Cedex 2, France}

\begin{document}

\maketitle

\begin{abstract} 
We propose an algebraic study of the simple graph
isomorphism problem. We define a Hopf algebra from an explicit
realization of its elements as formal power series. We show that these
series can be evaluated on graphs and count occurrences of subgraphs.
We establish a criterion for the isomorphism test of two simple graphs
by means of occurrence counting of subgraphs. This criterion is deduced
from algebraic relations between elements of our algebra.

\end{abstract}


\section{Introduction}

Let $G$ be a finite graph, drawn on a paper sheet. $G$ is represented
with dots for the vertices and segments between some pairs of
vertices. If $n$ counts the number of edges, if $a$ counts the number of
adjacent pairs of edges and $d$ the number of pairs of disjoint edges,
obviously, we have $a + d = \binom{n}{2}$. Now, we want to play
the same game with more complicated patterns and investigate how 
these counting functions are related. These
motivations are not new: similar approaches can be found for example 
in~\cite{MR1873218} or in~\cite{Mikko} but these papers
deal with finite cases where graphs have a fixed number of vertices.

We will take care of two difficult points for giving an algebraic sense to these counting
functions. The first one is that one can describe (for example on 
a computer) only labeled graphs, so we will have to deal with the
symmetries consisting in relabeling graphs. The second one is that we
do not want to depend on the size of graphs. Counting the number
of occurrences of a pattern inside a graph has a meaning whatever the
size of the graph. To achieve this goal, we will make sense of these
counting function by means of elements of a Hopf algebra $UGQSym$
(Unlabeled Graph Quasi Symmetric functions) which is close
to the combinatorial Hopf algebra $GTSym$ introduced by Novelli,
Thibon and Thiéry~\cite{MR2103196}.

We will first realize functions counting occurrences of subgraph, as
power series in an infinite number of variables. In Section~\ref{trois},
we investigate the Hopf structure of the algebra $UGQSym$. We show how this
algebra is connected to invariants of graphs in Section~\ref{quatre} and
give a sufficient criterion for two graphs to be isomorphic by means of
subgraph occurrences counting. Finally, we apply our results to the reconstruction
problem of finite graphs.

\section{Subgraph enumeration functions as formal power series}~\label{deux}

Let $\mathbb{A}$ be an infinite set of variables $x_{ij}$ indexed by pairs
$(i,j)$ of positive integers such that $i < j$. An adapted
presentation of this alphabet can be made as a triangle as follows:

\begin{displaymath}
  \begin{array}{ccc|c|cccc}
    1     &       &       &       &       &       &       & \\
    x_{12} & 2     &       &       &       &       &       & \\
    x_{13} & x_{23} & 3     &       &       &       &       & \\ \hline
    x_{14} & x_{24} & x_{34} & 4     &       &       &       & \\ \hline
    x_{15} & x_{25} & x_{35} & x_{45} & 5     &       &       & \\
    x_{16} & x_{26} & x_{36} & x_{46} & x_{56} & 6     &       & \\
    x_{17} & x_{27} & x_{37} & x_{47} & x_{57} & x_{67} & 7     &   \\
    \vdots &   &  & \vdots &  &  &  &  \ddots \\
  \end{array}
\end{displaymath}


\begin{definition}
Let $n$ be an integer and $\sigma$ be a permutation of the symmetric group
$\SG_n$. We define an action of $\SG_n$ over pairs $(i, j)$ (variables $x_{ij}$) such that
$1 \leqslant i < j \leqslant n$ as follows
\begin{displaymath}
  \sigma \cdot (i, j) := 
\left\{ \begin{array}{lr}
(\sigma(i), \sigma(j)) & \text{if }~ \sigma(i) < \sigma(j), \\
(\sigma(j), \sigma(i)) & \text{otherwise}.
\end{array}\right. \qquad
 \left( \sigma \cdot x_{ij} := 
\left\{ \begin{array}{lr}
x_{\sigma(i)\sigma(j)} & \text{if }~ \sigma(i) < \sigma(j), \\
x_{\sigma(j)\sigma(i)} & \text{otherwise}.
\end{array}\right.\right)
\end{displaymath}
\end{definition}

For finite graphs over $n$ vertices, this action corresponds exactly
to the relabeling action on graphs viewed on the lower triangular parts
of their adjacency matrices.

\begin{definition}
Let $n \geqslant 2$ be an integer and $G$ be a graph over $n$ vertices
without isolated vertex. We denote by $((a_{ij}))_{1 \leqslant i < j \leqslant n}$ 
the lower triangular part of the incidence matrix of $G$ when one has
chosen any labeling of the vertices of $G$. We thus define an
\emph{invariant power series} $\mathcal{M}_G$ as 
\begin{displaymath}
\mathcal{M}_G := \sum_{r_1 < \dots < r_n} \left( \sum_{\sigma \in Orb_{\SG_n}(G)} \left( \prod_{1 \leqslant i < j \leqslant n} x_{r_i r_j}^{a_{\sigma \cdot (i,j)}} \right) \right),
\end{displaymath}
where $Orb_{\SG_n}(G)$ is a set of permutations of $\SG_n$ required
to deploy the orbit of the graph $G$ under the relabeling action.
For the unique graph over $1$ node (loops are not allowed, this graph cannot contain any edge), we set
$\mathcal{M}_{{\newcommand{\nodea}{\node[draw, circle, scale=0.2, fill] (a) {$$};}
\begin{tikzpicture}[auto]
\matrix[column sep=.0cm, row sep=.1cm,ampersand replacement=\&]{
         \nodea \\ 
};
\end{tikzpicture}}} := 1$.
\end{definition}

Here are the first examples of functions $\mathcal{M}_G$.
\begin{displaymath}
\mathcal{M}_{{\newcommand{\nodea}{\node[draw,circle,scale=0.2,fill] (a) {$$}
;}\newcommand{\nodeb}{\node[draw,circle,scale=0.2,fill] (b) {$$}
;}\begin{tikzpicture}[auto]
\matrix[column sep=.0cm, row sep=.1cm,ampersand replacement=\&]{
                  \& \nodeb  \\ 
         \nodea   \&         \\ 
};
\path[thick] (a) edge (b);
\end{tikzpicture}}} = \sum_{0<i<j} x_{ij} = x_{12} + x_{13} + x_{23} + x_{14} + x_{24} + x_{34} + \dots
\end{displaymath}

\begin{displaymath}
\mathcal{M}_{{\newcommand{\nodea}{\node[draw,circle,scale=0.2,fill] (a) {$$};}
\newcommand{\nodeb}{\node[draw,circle,scale=0.2,fill] (b) {$$};}
\newcommand{\nodec}{\node[draw,circle,scale=0.2,fill] (c) {$$};}
\begin{tikzpicture}[auto]
\matrix[column sep=.0cm, row sep=.1cm,ampersand replacement=\&]{
                  \& \nodeb \&  \\ 
         \nodea   \&        \& \nodec \\ 
};
\path[thick] (a) edge (b)
             (b) edge (c);
\end{tikzpicture}}} = \sum_{0<i<j<k} x_{ij}x_{ik} + x_{ij}x_{jk} + x_{ik}x_{jk} = x_{12}x_{13} + x_{12}x_{23} + x_{13}x_{23} + \dots
\end{displaymath}

\begin{displaymath}
\mathcal{M}_{{\newcommand{\nodea}{\node[draw,circle,scale=0.2,fill] (a) {$$};}
\newcommand{\nodeb}{\node[draw,circle,scale=0.2,fill] (b) {$$};}
\newcommand{\nodec}{\node[draw,circle,scale=0.2,fill] (c) {$$};}
\newcommand{\noded}{\node[draw,circle,scale=0.2,fill] (d) {$$};}
\begin{tikzpicture}[auto]
\matrix[column sep=.0cm, row sep=.1cm,ampersand replacement=\&]{
                  \& \nodeb \&        \& \noded \\ 
         \nodea   \&        \& \nodec \& \\ 
};
\path[thick] (a) edge (b)
             (c) edge (d);
\end{tikzpicture}}} = \sum_{0<i<j<k<l} x_{ij}x_{kl} + x_{ik}x_{jl} + x_{il}x_{jk} = x_{12}x_{34} + x_{13}x_{24} + x_{14}x_{23} + \dots
\end{displaymath}

For each finite labeled graph $g$ with positive integers, we associate a
monomial $m(g)$ with $g$ as follows
\begin{displaymath}
  m(g) := \prod_{0 < i < j \atop{\text{ i and j are linked}}} x_{ij}
\end{displaymath}

\begin{proposition}
For $n \geqslant 2$ an integer and $G$ a graph over $n$ vertices
without isolated vertex, an alternative description of $\mathcal{M}_G$
is given by
\begin{displaymath}
\mathcal{M}_G = \sum_{g : \text{labeling of } G \atop{\text{ using positive integers}}} m(g).
\end{displaymath}
\begin{proof}
Monomials are in bijection with labeled graphs. The $\mathcal{M}_G$ are
infinite sums over some monomials and we can set a canonical
representative for each $\mathcal{M}_G$ using two commutative congruence relations
over the commutative monoid $(x_{ij} | x_{ij}^2 - x_{ij})^*$. As described in details in~\cite{MR3091039}, 
a canonical monomial is the maximum in its orbit under the relabeling action for the lexicographic order.
\begin{displaymath}
\begin{CD}
g_{support}  @>pack>> g_{1 \dots n} \\
@VcanonicalVV          @VVcanonicalV \\
G_{support}  @>pack>> G_{1 \dots n} \\
\end{CD}
\end{displaymath}
For example, we take the monomial (or labeled graph) $x_{25}x_{57}$.
\begin{displaymath}
\begin{CD}
x_{25}x_{57}  @>(2,5,7) \rightarrow (1,2,3)>pack> x_{12}x_{23} \\
@V(2,5,7) \rightarrow (5,2,7)VcanonicalV          @VcanonicalV(1,2,3) \rightarrow (2,1,3)V \\
x_{25}x_{27}  @>(2,5,7) \rightarrow (1,2,3)>pack> x_{12}x_{13} \\
\end{CD}
\end{displaymath}
\end{proof}
\end{proposition}

Let $H$ be a finite graph over $n$ nodes. Let $(a_{ij})_{1 \leqslant i < j \leqslant n}$ be 
the lower triangular part of its incidence matrix corresponding to any
labeling of $H$ with integers $\{1, \dots , n\}$. Now complete the
triangle $a_{ij}$ with an infinite number of $0$ for $a_{ij}$ where
$j > n$. The evaluation of $\mathcal{M}_G$ over the infinite sequence
of boolean associated with $H$ counts the number of embeddings of $G$
inside $H$.
\begin{displaymath}
  \mathcal{M}_G(H) := \mathcal{M}_G((a_{ij})) = \#\{ G \hookrightarrow H \} \in \mathbb{N}
\end{displaymath}
\begin{displaymath}
  \mathcal{M}_G(G) = 1 \qquad \forall G \text{ graph without isolated vertex}
\end{displaymath}

The main motivation justifying the definition of $\mathcal{M}_G$ stand in the following result.

\begin{theorem}
$\mathcal{M}_G$ can be seen as a maps, defined on every finite graph, counting
occurrences of subgraphs. For example, the function 
$\mathcal{M}_{{\newcommand{\nodea}{\node[draw,circle,scale=0.2,fill] (a) {$$};}
\newcommand{\nodeb}{\node[draw,circle,scale=0.2,fill] (b) {$$};}
\newcommand{\nodec}{\node[draw,circle,scale=0.2,fill] (c) {$$};}
\begin{tikzpicture}[auto]
\matrix[column sep=.0cm, row sep=.1cm,ampersand replacement=\&]{
                  \& \nodeb \&        \\ 
         \nodea   \&        \& \nodec \\ 
};
\path[thick] (a) edge (b) 
             (b) edge (c)
             (c) edge (a);
\end{tikzpicture}}}$ counts the number of triangles inside a graph.
\end{theorem}

$\mathcal{M}_{{\newcommand{\nodea}{\node[draw,circle,scale=0.2,fill] (a) {$$};}
\newcommand{\nodeb}{\node[draw,circle,scale=0.2,fill] (b) {$$};}
\newcommand{\nodec}{\node[draw,circle,scale=0.2,fill] (c) {$$};}
\begin{tikzpicture}[auto]
\matrix[column sep=.0cm, row sep=.1cm,ampersand replacement=\&]{
                  \& \nodeb \&        \\ 
         \nodea   \&        \& \nodec \\ 
};
\path[thick] (a) edge (b) 
             (b) edge (c)
             (c) edge (a);
\end{tikzpicture}}} = x_{12}x_{13}x_{23} + x_{12}x_{14}x_{24} + x_{13}x_{14}x_{34} + x_{23}x_{24}x_{34} + x_{12}x_{15}x_{25} +
x_{13}x_{15}x_{35} + x_{23}x_{25}x_{35} + x_{14}x_{15}x_{45} + x_{24}x_{25}x_{45} + x_{34}x_{35}x_{45} + \dots$

If $H$ is the following graph over $5$ nodes, we give it labels with integers from $1$ to $5$
\vspace{-0.3cm}
\begin{displaymath}
    \newcommand{\nodea}{\node[draw,circle,scale=1] (a) {$\textcolor{white}{1}$};}
    \newcommand{\nodeb}{\node[draw,circle,scale=1] (b) {$\textcolor{white}{2}$};}
    \newcommand{\nodec}{\node[draw,circle,scale=1] (c) {$\textcolor{white}{3}$};}
    \newcommand{\noded}{\node[draw,circle,scale=1] (d) {$\textcolor{white}{4}$};}
    \newcommand{\nodee}{\node[draw,circle,scale=1] (e) {$\textcolor{white}{5}$};}
    \begin{tikzpicture}[auto]
      \matrix[column sep=.1cm, row sep=.3cm,ampersand replacement=\&]{
               \&        \& \nodeb \&        \& \\ 
        \nodea \&        \&        \&        \& \nodec \\ 
               \& \noded \&        \& \nodee \& \\
      };
      \path[thick] (a) edge (b) 
      (b) edge (c)
      (c) edge (a)
      (a) edge (d)
      (b) edge (d)
      (a) edge (e)
      (c) edge (e);
    \end{tikzpicture} \qquad
    \renewcommand{\nodea}{\node[draw,circle,scale=1] (a) {$1$};}
    \renewcommand{\nodeb}{\node[draw,circle,scale=1] (b) {$2$};}
    \renewcommand{\nodec}{\node[draw,circle,scale=1] (c) {$3$};}
    \renewcommand{\noded}{\node[draw,circle,scale=1] (d) {$4$};}
    \renewcommand{\nodee}{\node[draw,circle,scale=1] (e) {$5$};}
    \begin{tikzpicture}[auto]
      \matrix[column sep=.1cm, row sep=.3cm,ampersand replacement=\&]{
               \&        \& \nodeb \&        \& \\ 
        \nodea \&        \&        \&        \& \nodec \\ 
               \& \noded \&        \& \nodee \& \\
      };
      \path[thick] (a) edge (b) 
      (b) edge (c)
      (c) edge (a)
      (a) edge (d)
      (b) edge (d)
      (a) edge (e)
      (c) edge (e);
    \end{tikzpicture}
\end{displaymath}
\vspace{-0.7cm}
\begin{displaymath}
   \text{incidence matrix}(H) : \qquad \quad
    \begin{array}{c|ccccc}
      & \mathbf{1} & \mathbf{2} & \mathbf{3} & \mathbf{4} & \mathbf{5} \\ \hline
      \mathbf{1} & \cdot \\
      \mathbf{2} & 1 & \cdot \\
      \mathbf{3} & 1 & 1 & \cdot \\
      \mathbf{4} & 1 & 1 & 0 & \cdot \\
      \mathbf{5} & 1 & 0 & 1 & 0 & \cdot \\
    \end{array}
\end{displaymath}

$\mathcal{M}_{{\newcommand{\nodea}{\node[draw,circle,scale=0.2,fill] (a) {$$};}
\newcommand{\nodeb}{\node[draw,circle,scale=0.2,fill] (b) {$$};}
\newcommand{\nodec}{\node[draw,circle,scale=0.2,fill] (c) {$$};}
\begin{tikzpicture}[auto]
\matrix[column sep=.0cm, row sep=.1cm,ampersand replacement=\&]{
                  \& \nodeb \&        \\ 
         \nodea   \&        \& \nodec \\ 
};
\path[thick] (a) edge (b) 
             (b) edge (c)
             (c) edge (a);
\end{tikzpicture}}}(H) = 1 \cdot 1 \cdot 1 + 1 \cdot 1 \cdot 1 + 
1 \cdot 1 \cdot 0 + 1 \cdot 1 \cdot 0 + 1 \cdot 1 \cdot 0 +
1 \cdot 1 \cdot 1 + 1 \cdot 0 \cdot 1 + 1 \cdot 1 \cdot 0 + 
1 \cdot 0 \cdot 0 + 0 \cdot 1 \cdot 0 + \dots$

$\mathcal{M}_{{\newcommand{\nodea}{\node[draw,circle,scale=0.2,fill] (a) {$$};}
\newcommand{\nodeb}{\node[draw,circle,scale=0.2,fill] (b) {$$};}
\newcommand{\nodec}{\node[draw,circle,scale=0.2,fill] (c) {$$};}
\begin{tikzpicture}[auto]
\matrix[column sep=.0cm, row sep=.1cm,ampersand replacement=\&]{
                  \& \nodeb \&        \\ 
         \nodea   \&        \& \nodec \\ 
};
\path[thick] (a) edge (b) 
             (b) edge (c)
             (c) edge (a);
\end{tikzpicture}}}(H) = 3$ as the only three monomials contributing 
to $H$ are $x_{12}x_{13}x_{23}$, $x_{12}x_{14}x_{24}$ and $x_{13}x_{15}x_{35}$. 
This number of triangles is independent of the way chosen for labeling
the graph since the functions $\mathcal{M}_G$ are invariant under the
relabeling action.

\section{The algebra $UGQSym$}~\label{trois}

We denote by $UGQSym$ the subspace of $\mathbb{K}[x_{ij} | x_{ij}^2=x_{ij}, i<j]$ 
generated by the $\mathcal{M}_G$. We will call this space the 
\emph{Unlabeled Graph Quasi Symmetric functions}.

As we focus on simple graphs (unoriented, simple-edged), the
variables, which model edges, can take only two values: $0$ if the
edge is not present and $1$ otherwise. These two numbers are the roots
of the polynomial $x^2 - x$ and choosing formal variables $x_{ij}$
satisfying $x_{ij}^2 = x_{ij}$ allow us to describe things more
combinatorically.

\begin{theorem}
$UGQSym$ is a subalgebra of $\mathbb{K}[x_{ij} | i<j : x_{ij}^2 = x_{ij}]$. 
Precisely, there exist non negative integers $c_{G', G''}^G$ such that
\begin{displaymath}
\mathcal{M}_{G'} \cdot \mathcal{M}_{G''} = \sum_G c_{G', G''}^G \mathcal{M}_G.
\end{displaymath}
The $c_{G', G''}^G$ counts precisely the number of ways of embedding
simultaneously $G'$ and $G''$ onto $G$.
\end{theorem}

Here are some examples of products obtained with Sage~\cite{sage} using 
tools for computing canonicals under relabeling action~\cite{MR3091039}.

\begin{displaymath}
\mathcal{M}_{{\newcommand{\nodea}{\node[draw,circle,scale=0.2,fill] (a) {$$};}
\newcommand{\nodeb}{\node[draw,circle,scale=0.2,fill] (b) {$$};}
\begin{tikzpicture}[auto]
\matrix[column sep=.0cm, row sep=.1cm,ampersand replacement=\&]{
                  \& \nodeb \\ 
         \nodea   \&        \\ 
};
\path[thick] (a) edge (b);
\end{tikzpicture}}} \cdot
\mathcal{M}_{{\newcommand{\nodea}{\node[draw,circle,scale=0.2,fill] (a) {$$};}
\newcommand{\nodeb}{\node[draw,circle,scale=0.2,fill] (b) {$$};}
\begin{tikzpicture}[auto]
\matrix[column sep=.0cm, row sep=.1cm,ampersand replacement=\&]{
                  \& \nodeb \\ 
         \nodea   \&        \\ 
};
\path[thick] (a) edge (b);
\end{tikzpicture}}} = 
\mathcal{M}_{{\newcommand{\nodea}{\node[draw,circle,scale=0.2,fill] (a) {$$};}
\newcommand{\nodeb}{\node[draw,circle,scale=0.2,fill] (b) {$$};}
\begin{tikzpicture}[auto]
\matrix[column sep=.0cm, row sep=.1cm,ampersand replacement=\&]{
                  \& \nodeb \\ 
         \nodea   \&        \\ 
};
\path[thick] (a) edge (b);
\end{tikzpicture}}} +
2 \mathcal{M}_{{\newcommand{\nodea}{\node[draw,circle,scale=0.2,fill] (a) {$$};}
\newcommand{\nodeb}{\node[draw,circle,scale=0.2,fill] (b) {$$};}
\newcommand{\nodec}{\node[draw,circle,scale=0.2,fill] (c) {$$};}
\begin{tikzpicture}[auto]
\matrix[column sep=.0cm, row sep=.1cm,ampersand replacement=\&]{
                  \& \nodeb \&        \\ 
         \nodea   \&        \& \nodec \\ 
};
\path[thick] (a) edge (b)
             (b) edge (c);
\end{tikzpicture}}} +
2 \mathcal{M}_{{\newcommand{\nodea}{\node[draw,circle,scale=0.2,fill] (a) {$$};}
\newcommand{\nodeb}{\node[draw,circle,scale=0.2,fill] (b) {$$};}
\newcommand{\nodec}{\node[draw,circle,scale=0.2,fill] (c) {$$};}
\newcommand{\noded}{\node[draw,circle,scale=0.2,fill] (d) {$$};}
\begin{tikzpicture}[auto]
\matrix[column sep=.0cm, row sep=.1cm,ampersand replacement=\&]{
                  \& \nodeb \&        \& \noded \\ 
         \nodea   \&        \& \nodec \&        \\ 
};
\path[thick] (a) edge (b)
             (c) edge (d);
\end{tikzpicture}}}
\end{displaymath}

\begin{equation}~\label{prod1}
\mathcal{M}_{{\newcommand{\nodea}{\node[draw,circle,scale=0.2,fill] (a) {$$};}
\newcommand{\nodeb}{\node[draw,circle,scale=0.2,fill] (b) {$$};}
\newcommand{\nodec}{\node[draw,circle,scale=0.2,fill] (c) {$$};}
\begin{tikzpicture}[auto]
\matrix[column sep=.0cm, row sep=.1cm,ampersand replacement=\&]{
                  \& \nodeb \&  \\ 
         \nodea   \&        \& \nodec \\ 
};
\path[thick] (a) edge (b) 
             (b) edge (c);
\end{tikzpicture}}} \cdot
\mathcal{M}_{{\newcommand{\nodea}{\node[draw,circle,scale=0.2,fill] (a) {$$};}
\newcommand{\nodeb}{\node[draw,circle,scale=0.2,fill] (b) {$$};}
\begin{tikzpicture}[auto]
\matrix[column sep=.0cm, row sep=.1cm,ampersand replacement=\&]{
                  \& \nodeb \\ 
         \nodea   \&        \\ 
};
\path[thick] (a) edge (b);
\end{tikzpicture}}} = 
2 \mathcal{M}_{{\newcommand{\nodea}{\node[draw,circle,scale=0.2,fill] (a) {$$};}
\newcommand{\nodeb}{\node[draw,circle,scale=0.2,fill] (b) {$$};}
\newcommand{\nodec}{\node[draw,circle,scale=0.2,fill] (c) {$$};}
\begin{tikzpicture}[auto]
\matrix[column sep=.0cm, row sep=.1cm,ampersand replacement=\&]{
                  \& \nodeb \&  \\ 
         \nodea   \&        \& \nodec \\ 
};
\path[thick] (a) edge (b) 
             (b) edge (c);
\end{tikzpicture}}} +
3 \mathcal{M}_{{\newcommand{\nodea}{\node[draw,circle,scale=0.2,fill] (a) {$$};}
\newcommand{\nodeb}{\node[draw,circle,scale=0.2,fill] (b) {$$};}
\newcommand{\nodec}{\node[draw,circle,scale=0.2,fill] (c) {$$};}
\begin{tikzpicture}[auto]
\matrix[column sep=.0cm, row sep=.1cm,ampersand replacement=\&]{
                  \& \nodeb \& \\ 
         \nodea   \&        \& \nodec \\ 
};
\path[thick] (a) edge (b)
             (b) edge (c)
             (c) edge (a);
\end{tikzpicture}}} +
3 \mathcal{M}_{{\newcommand{\nodea}{\node[draw,circle,scale=0.2,fill] (a) {$$};}
\newcommand{\nodeb}{\node[draw,circle,scale=0.2,fill] (b) {$$};}
\newcommand{\nodec}{\node[draw,circle,scale=0.2,fill] (c) {$$};}
\newcommand{\noded}{\node[draw,circle,scale=0.2,fill] (d) {$$};}
\begin{tikzpicture}[auto]
\matrix[column sep=.03cm, row sep=.04cm,ampersand replacement=\&]{
                  \& \noded \& \\
                  \& \nodea \& \\ 
         \nodeb   \&        \& \nodec \\ 
};
\path[thick] (a) edge (b)
             (a) edge (c)
             (a) edge (d);
\end{tikzpicture}}} +
2 \mathcal{M}_{{\newcommand{\nodea}{\node[draw,circle,scale=0.2,fill] (a) {$$};}
\newcommand{\nodeb}{\node[draw,circle,scale=0.2,fill] (b) {$$};}
\newcommand{\nodec}{\node[draw,circle,scale=0.2,fill] (c) {$$};}
\newcommand{\noded}{\node[draw,circle,scale=0.2,fill] (d) {$$};}
\begin{tikzpicture}[auto]
\matrix[column sep=.0cm, row sep=.1cm,ampersand replacement=\&]{
                  \& \nodeb \&        \& \noded \\ 
         \nodea   \&        \& \nodec \& \\ 
};
\path[thick] (a) edge (b) 
             (b) edge (c)
             (c) edge (d);
\end{tikzpicture}}} + 
\mathcal{M}_{{\newcommand{\nodea}{\node[draw,circle,scale=0.2,fill] (a) {$$};}
\newcommand{\nodeb}{\node[draw,circle,scale=0.2,fill] (b) {$$};}
\newcommand{\nodec}{\node[draw,circle,scale=0.2,fill] (c) {$$};}
\newcommand{\noded}{\node[draw,circle,scale=0.2,fill] (d) {$$};}
\newcommand{\nodee}{\node[draw,circle,scale=0.2,fill] (e) {$$};}
\begin{tikzpicture}[auto]
\matrix[column sep=.0cm, row sep=.1cm,ampersand replacement=\&]{
                  \& \nodeb \&        \& \noded \& \\ 
         \nodea   \&        \& \nodec \&        \& \nodee \\ 
};
\path[thick] (a) edge (b) 
             (b) edge (c)
             (d) edge (e);
\end{tikzpicture}}}
\end{equation}

\begin{displaymath}
\mathcal{M}_{{\newcommand{\nodea}{\node[draw,circle,scale=0.2,fill] (a) {$$};}
\newcommand{\nodeb}{\node[draw,circle,scale=0.2,fill] (b) {$$};}
\newcommand{\nodec}{\node[draw,circle,scale=0.2,fill] (c) {$$};}
\newcommand{\noded}{\node[draw,circle,scale=0.2,fill] (d) {$$};}
\begin{tikzpicture}[auto]
\matrix[column sep=.0cm, row sep=.1cm,ampersand replacement=\&]{
                  \& \nodeb \&        \& \noded \\ 
         \nodea   \&        \& \nodec \&        \\ 
};
\path[thick] (a) edge (b) 
             (c) edge (d);
\end{tikzpicture}}} \cdot 
\mathcal{M}_{{\newcommand{\nodea}{\node[draw,circle,scale=0.2,fill] (a) {$$};}
\newcommand{\nodeb}{\node[draw,circle,scale=0.2,fill] (b) {$$};}
\begin{tikzpicture}[auto]
\matrix[column sep=.0cm, row sep=.1cm,ampersand replacement=\&]{
                  \& \nodeb \\ 
         \nodea   \&        \\ 
};
\path[thick] (a) edge (b);
\end{tikzpicture}}} = 
2 \mathcal{M}_{{\newcommand{\nodea}{\node[draw,circle,scale=0.2,fill] (a) {$$};}
\newcommand{\nodeb}{\node[draw,circle,scale=0.2,fill] (b) {$$};}
\newcommand{\nodec}{\node[draw,circle,scale=0.2,fill] (c) {$$};}
\newcommand{\noded}{\node[draw,circle,scale=0.2,fill] (d) {$$};}
\begin{tikzpicture}[auto]
\matrix[column sep=.0cm, row sep=.1cm,ampersand replacement=\&]{
                  \& \nodeb \&        \& \noded \\ 
         \nodea   \&        \& \nodec \&        \\ 
};
\path[thick] (a) edge (b) 
             (c) edge (d);
\end{tikzpicture}}} +
\mathcal{M}_{{\newcommand{\nodea}{\node[draw,circle,scale=0.2,fill] (a) {$$};}
\newcommand{\nodeb}{\node[draw,circle,scale=0.2,fill] (b) {$$};}
\newcommand{\nodec}{\node[draw,circle,scale=0.2,fill] (c) {$$};}
\newcommand{\noded}{\node[draw,circle,scale=0.2,fill] (d) {$$};}
\begin{tikzpicture}[auto]
\matrix[column sep=.0cm, row sep=.1cm,ampersand replacement=\&]{
                  \& \nodeb \&        \& \noded \\ 
         \nodea   \&        \& \nodec \&        \\ 
};
\path[thick] (a) edge (b) 
             (b) edge (c)
             (c) edge (d);
\end{tikzpicture}}} +
2 \mathcal{M}_{{\newcommand{\nodea}{\node[draw,circle,scale=0.2,fill] (a) {$$};}
\newcommand{\nodeb}{\node[draw,circle,scale=0.2,fill] (b) {$$};}
\newcommand{\nodec}{\node[draw,circle,scale=0.2,fill] (c) {$$};}
\newcommand{\noded}{\node[draw,circle,scale=0.2,fill] (d) {$$};}
\newcommand{\nodee}{\node[draw,circle,scale=0.2,fill] (e) {$$};}
\begin{tikzpicture}[auto]
\matrix[column sep=.0cm, row sep=.1cm,ampersand replacement=\&]{
                  \& \nodeb \&        \& \noded \&        \\ 
         \nodea   \&        \& \nodec \&        \& \nodee \\ 
};
\path[thick] (a) edge (b) 
             (b) edge (c)
             (d) edge (e);
\end{tikzpicture}}} +
3 \mathcal{M}_{{\newcommand{\nodea}{\node[draw,circle,scale=0.2,fill] (a) {$$};}
\newcommand{\nodeb}{\node[draw,circle,scale=0.2,fill] (b) {$$};}
\newcommand{\nodec}{\node[draw,circle,scale=0.2,fill] (c) {$$};}
\newcommand{\noded}{\node[draw,circle,scale=0.2,fill] (d) {$$};}
\newcommand{\nodee}{\node[draw,circle,scale=0.2,fill] (e) {$$};}
\newcommand{\nodef}{\node[draw,circle,scale=0.2,fill] (f) {$$};}
\begin{tikzpicture}[auto]
\matrix[column sep=.0cm, row sep=.1cm,ampersand replacement=\&]{
                  \& \nodeb \&        \& \noded \&        \& \nodef \\ 
         \nodea   \&        \& \nodec \&        \& \nodee \&        \\ 
};
\path[thick] (a) edge (b) 
             (c) edge (d)
             (e) edge (f);
\end{tikzpicture}}}
\end{displaymath}

\begin{displaymath}
\mathcal{M}_{{\newcommand{\nodea}{\node[draw,circle,scale=0.2,fill] (a) {$$};}
\newcommand{\nodeb}{\node[draw,circle,scale=0.2,fill] (b) {$$};}
\newcommand{\nodec}{\node[draw,circle,scale=0.2,fill] (c) {$$};}
\begin{tikzpicture}[auto]
\matrix[column sep=.0cm, row sep=.1cm,ampersand replacement=\&]{
                  \& \nodeb \&        \\ 
         \nodea   \&        \& \nodec \\ 
};
\path[thick] (a) edge (b) 
             (b) edge (c)
             (c) edge (a);
\end{tikzpicture}}} \cdot  
\mathcal{M}_{{\newcommand{\nodea}{\node[draw,circle,scale=0.2,fill] (a) {$$};}
\newcommand{\nodeb}{\node[draw,circle,scale=0.2,fill] (b) {$$};}
\newcommand{\nodec}{\node[draw,circle,scale=0.2,fill] (c) {$$};}
\begin{tikzpicture}[auto]
\matrix[column sep=.0cm, row sep=.1cm,ampersand replacement=\&]{
                  \& \nodeb \&        \\ 
         \nodea   \&        \& \nodec \\ 
};
\path[thick] (a) edge (b) 
             (b) edge (c)
             (c) edge (a);
\end{tikzpicture}}} = 
\mathcal{M}_{{\newcommand{\nodea}{\node[draw,circle,scale=0.2,fill] (a) {$$};}
\newcommand{\nodeb}{\node[draw,circle,scale=0.2,fill] (b) {$$};}
\newcommand{\nodec}{\node[draw,circle,scale=0.2,fill] (c) {$$};}
\begin{tikzpicture}[auto]
\matrix[column sep=.0cm, row sep=.1cm,ampersand replacement=\&]{
                  \& \nodeb \&        \\ 
         \nodea   \&        \& \nodec \\ 
};
\path[thick] (a) edge (b) 
             (b) edge (c)
             (c) edge (a);
\end{tikzpicture}}} +
2 \mathcal{M}_{{\newcommand{\nodea}{\node[draw,circle,scale=0.2,fill] (a) {$$};}
\newcommand{\nodeb}{\node[draw,circle,scale=0.2,fill] (b) {$$};}
\newcommand{\nodec}{\node[draw,circle,scale=0.2,fill] (c) {$$};}
\newcommand{\noded}{\node[draw,circle,scale=0.2,fill] (d) {$$};}
\begin{tikzpicture}[auto]
\matrix[column sep=.0cm, row sep=.1cm,ampersand replacement=\&]{
                  \& \nodeb \&        \& \noded \\ 
         \nodea   \&        \& \nodec \&        \\ 
};
\path[thick] (a) edge (b) 
             (b) edge (c)
             (c) edge (a)
             (b) edge (d)
             (c) edge (d);
\end{tikzpicture}}} +
2 \mathcal{M}_{{\newcommand{\nodea}{\node[draw,circle,scale=0.2,fill] (a) {$$};}
\newcommand{\nodeb}{\node[draw,circle,scale=0.2,fill] (b) {$$};}
\newcommand{\nodec}{\node[draw,circle,scale=0.2,fill] (c) {$$};}
\newcommand{\noded}{\node[draw,circle,scale=0.2,fill] (d) {$$};}
\newcommand{\nodee}{\node[draw,circle,scale=0.2,fill] (e) {$$};}
\begin{tikzpicture}[auto]
\matrix[column sep=.0cm, row sep=.1cm,ampersand replacement=\&]{
                  \& \nodeb \&        \& \noded \&        \\ 
         \nodea   \&        \& \nodec \&        \& \nodee \\ 
};
\path[thick] (a) edge (b) 
             (b) edge (c)
             (c) edge (a)
             (c) edge (d)
             (c) edge (e)
             (e) edge (d);
\end{tikzpicture}}} +
2 \mathcal{M}_{{\newcommand{\nodea}{\node[draw,circle,scale=0.2,fill] (a) {$$};}
\newcommand{\nodeb}{\node[draw,circle,scale=0.2,fill] (b) {$$};}
\newcommand{\nodec}{\node[draw,circle,scale=0.2,fill] (c) {$$};}
\newcommand{\noded}{\node[draw,circle,scale=0.2,fill] (d) {$$};}
\newcommand{\nodee}{\node[draw,circle,scale=0.2,fill] (e) {$$};}
\newcommand{\nodef}{\node[draw,circle,scale=0.2,fill] (f) {$$};}
\begin{tikzpicture}[auto]
\matrix[column sep=.0cm, row sep=.1cm,ampersand replacement=\&]{
                  \& \nodeb \&        \& \noded \&        \& \nodef \\ 
         \nodea   \&        \& \nodec \&        \& \nodee \&        \\ 
};
\path[thick] (a) edge (b) 
             (b) edge (c)
             (c) edge (a)
             (d) edge (e)
             (f) edge (e)
             (f) edge (d);
\end{tikzpicture}}}
\end{displaymath}

As an exercise, the reader can check that the definition of functions
$\mathcal{M}_G$ realizes the simple combinatorial fact: inside a graph,
pairs of edges are adjacent or disjoint (use $x_{ij}^2 = x_{ij}$).
\begin{displaymath}
\binom{\mathcal{M}_{{\newcommand{\nodea}{\node[draw,circle,scale=0.2,fill] (a) {$$}
;}\newcommand{\nodeb}{\node[draw,circle,scale=0.2,fill] (b) {$$}
;}\begin{tikzpicture}[auto]
\matrix[column sep=.0cm, row sep=.1cm,ampersand replacement=\&]{
                  \& \nodeb  \\ 
         \nodea   \&         \\ 
};
\path[thick] (a) edge (b);
\end{tikzpicture}}}}{2} :=
\frac{{(\mathcal{M}_{{\newcommand{\nodea}{\node[draw,circle,scale=0.2,fill] (a) {$$}
;}\newcommand{\nodeb}{\node[draw,circle,scale=0.2,fill] (b) {$$}
;}\begin{tikzpicture}[auto]
\matrix[column sep=.0cm, row sep=.1cm,ampersand replacement=\&]{
                  \& \nodeb  \\ 
         \nodea   \&         \\ 
};
\path[thick] (a) edge (b);
\end{tikzpicture}}})}
{(\mathcal{M}_{{\newcommand{\nodea}{\node[draw,circle,scale=0.2,fill] (a) {$$}
;}\newcommand{\nodeb}{\node[draw,circle,scale=0.2,fill] (b) {$$}
;}\begin{tikzpicture}[auto]
\matrix[column sep=.0cm, row sep=.1cm,ampersand replacement=\&]{
                  \& \nodeb  \\ 
         \nodea   \&         \\ 
};
\path[thick] (a) edge (b);
\end{tikzpicture}}} - 1)}}{2!} = 
\mathcal{M}_{{\newcommand{\nodea}{\node[draw,circle,scale=0.2,fill] (a) {$$};}
\newcommand{\nodeb}{\node[draw,circle,scale=0.2,fill] (b) {$$};}
\newcommand{\nodec}{\node[draw,circle,scale=0.2,fill] (c) {$$};}
\begin{tikzpicture}[auto]
\matrix[column sep=.0cm, row sep=.1cm,ampersand replacement=\&]{
                  \& \nodeb \&  \\ 
         \nodea   \&        \& \nodec \\ 
};
\path[thick] (a) edge (b) 
             (b) edge (c);
\end{tikzpicture}}} + 
\mathcal{M}_{{\newcommand{\nodea}{\node[draw,circle,scale=0.2,fill] (a) {$$};}
\newcommand{\nodeb}{\node[draw,circle,scale=0.2,fill] (b) {$$};}
\newcommand{\nodec}{\node[draw,circle,scale=0.2,fill] (c) {$$};}
\newcommand{\noded}{\node[draw,circle,scale=0.2,fill] (d) {$$};}
\begin{tikzpicture}[auto]
\matrix[column sep=.0cm, row sep=.1cm,ampersand replacement=\&]{
                  \& \nodeb \&        \& \noded \\ 
         \nodea   \&        \& \nodec \&   \\ 
};
\path[thick] (a) edge (b) 
             (c) edge (d);
\end{tikzpicture}}}
\end{displaymath}

More generally, it is possible to show that for any non negative integer $n$, we have
\begin{displaymath}
\binom{\mathcal{M}_{{\newcommand{\nodea}{\node[draw,circle,scale=0.2,fill] (a) {$$}
;}\newcommand{\nodeb}{\node[draw,circle,scale=0.2,fill] (b) {$$}
;}\begin{tikzpicture}[auto]
\matrix[column sep=.0cm, row sep=.1cm,ampersand replacement=\&]{
                  \& \nodeb  \\ 
         \nodea   \&         \\ 
};
\path[thick] (a) edge (b);
\end{tikzpicture}}}}{n} :=
\frac{{(\mathcal{M}_{{\newcommand{\nodea}{\node[draw,circle,scale=0.2,fill] (a) {$$}
;}\newcommand{\nodeb}{\node[draw,circle,scale=0.2,fill] (b) {$$}
;}\begin{tikzpicture}[auto]
\matrix[column sep=.0cm, row sep=.1cm,ampersand replacement=\&]{
                  \& \nodeb  \\ 
         \nodea   \&         \\ 
};
\path[thick] (a) edge (b);
\end{tikzpicture}}})}
{(\mathcal{M}_{{\newcommand{\nodea}{\node[draw,circle,scale=0.2,fill] (a) {$$}
;}\newcommand{\nodeb}{\node[draw,circle,scale=0.2,fill] (b) {$$}
;}\begin{tikzpicture}[auto]
\matrix[column sep=.0cm, row sep=.1cm,ampersand replacement=\&]{
                  \& \nodeb  \\ 
         \nodea   \&         \\ 
};
\path[thick] (a) edge (b);
\end{tikzpicture}}} - 1)} \dots 
{(\mathcal{M}_{{\newcommand{\nodea}{\node[draw,circle,scale=0.2,fill] (a) {$$}
;}\newcommand{\nodeb}{\node[draw,circle,scale=0.2,fill] (b) {$$}
;}\begin{tikzpicture}[auto]
\matrix[column sep=.0cm, row sep=.1cm,ampersand replacement=\&]{
                  \& \nodeb  \\ 
         \nodea   \&         \\ 
};
\path[thick] (a) edge (b);
\end{tikzpicture}}} - (n-1))}}{n!}
= \sum_{G \text{ has $n$ edges}} \mathcal{M}_G
\end{displaymath}

\begin{proposition}
The set of functions $\{\mathcal{M}_G\}$ for a graph $G$ without
isolated vertex forms a linear basis of the algebra $UGQSym$.

\begin{proof}
By definition, the family generates $UGQSym$. Defining leading terms as
canonical labeled graphs over packed support, we immediately see that
elements of this family are linearly independent.
\end{proof}
\end{proposition}

The products of functions $\{\mathcal{M}_G\}$ establish algebraic
relations between functions counting occurrences of subgraphs. Finally,
we can show that only connected patterns are important.

\begin{proposition}
The functions $\{\mathcal{M}_G\}$ for $G$ connected graph
generates the algebra $UGQSym$.

\begin{proof}
By induction on the number of connected components. Any product 
$\mathcal{M}_{G_1} \times \mathcal{M}_{G_2}$ with respectively 
$n_1$ and $n_2$ connected components in $G_1$ and $G_2$ takes the form
\begin{displaymath}
\mathcal{M}_{G_1} \times \mathcal{M}_{G_2} = \mathcal{M}_{G_1 \sqcup G_2} + 
\sum_{G \text{ has } n < n_1 + n_2\atop{\text{connected components}}} \mathcal{M}_{G}
\end{displaymath}
\end{proof}
\end{proposition}

Filtered by number of nodes, the dimensions of "pseudo" homogeneous
components is counted by $graph\_node(n)$, the number of graphs over
$n$ nodes, Sequence A000088 of the OEIS~\cite{Sloane}. We
recall the first values $1, 1, 2, 4, 11, 34, 156, 1044, \dots$. These
graphs can be seen as monomials over connected graphs, for that, we
take connected components as variables and exponents are respectively
the multiplicities of isomorphic connected components appearing in the
graph. Setting $connected\_node(n)$ as the number of connected graphs
over $n$ nodes beginning by $1, 1, 1, 2, 6, 21, 112, 853, \dots$ 
also sequence A001349 of the OEIS~\cite{Sloane}. We
have the following relation with power series:
\begin{displaymath}
\prod_{n > 0} \frac{1}{(1 - q^n)^{connected\_node(n)}} = \sum_{n \geqslant 0} graph\_node(n) q^n .
\end{displaymath}
Therefore, the $\{\mathcal{M}_G\}$ for $G$ connected graphs are
algebraically independent.

We can play the same game as before, counting, this time, graphs by
number of edges. The number $graph\_edge(n)$ of graphs with $n$ edges (and without isolated vertex)
is A000664 of the OEIS and begins with $1, 1, 2, 5, 11, 26, 68, 177, 497$ 
and the number of connected graphs $connected\_edge(n)$ with $n$
edges is A002905 of the OEIS: $1, 1, 1, 3, 5, 12, 30, 79, 227, 710, \dots$.
Although the OEIS does not mention that the first one is the Euler
transform of the second, with the same argument as before, we also
have:
\begin{displaymath}
\prod_{n > 0} \frac{1}{(1 - q^n)^{connected\_edge(n)}} = \sum_{n \geqslant 0} graph\_edge(n) q^n .
\end{displaymath}

We can realize a coproduct by means of \emph{doubling alphabet trick} of~\cite{MR1935570}. 
Noting $\mathbb{A}$ our infinite triangle
of variables $x_{ij}$ with $0<i<j$, we introduce another copy of this
alphabet, $\mathbb{B}$ composed of variables $x_{i'j'}$ with $i'$ and
$j'$ satisfying the same constraints as $i$ and $j$. We define 
$\mathbb{A} \oplus \mathbb{B}$ like $\mathbb{A} \sqcup \mathbb{B}$ since we want
to forbid variables of the form $x_{ij'}$ or $x_{j'i}$.

\begin{displaymath}
  \begin{array}{rcl}
  \mathcal{M}_G(\mathbb{A} \oplus \mathbb{B}) & := & \displaystyle\sum_{g : \text{labeling of } G \atop{\text{ using normal and prime letters}}} m(g) \\
  & = & \displaystyle\sum_{g_1, g_2 : \text{labeling of } G_1 \sqcup G_2 = G \atop{\text{ $G_1$ labeled with normal letters, $G_2$ labeled with prime letters}}} m(g_1 \sqcup g_2)\\
  & = & \displaystyle\sum_{G_1 \sqcup G_2 = G}
(\displaystyle\sum_{g_1 : \text{labeling of } G_1 \atop{\text{ $G_1$ labeled with normal letters}}} m(g_1)) \cdot
(\displaystyle\sum_{g_2 : \text{labeling of } G_2 \atop{\text{ $G_2$ labeled with prime letters}}} m(g_2)) \\
  \mathcal{M}_G(\mathbb{A} \oplus \mathbb{B}) & = & \displaystyle\sum_{G_1 \sqcup G_2 = G} \mathcal{M}_{G_1}(\mathbb{A}) \otimes \mathbb{M}_{G_2}(\mathbb{B}) \\
  \end{array}
\end{displaymath}

\begin{definition}
Let $G$ be a finite graph without isolated vertex. We define a coproduct
$\Delta$ on basis elements $\mathcal{M}_G$ of $UGQSym$ as follows:
\begin{displaymath}
\Delta(\mathcal{M}_G) := \sum_{G' \sqcup G'' = G} M_{G'} \otimes M_{G''}
\end{displaymath}
The sum runs over \underline{ordered} pairs $(G', G'')$ of graphs without isolated
vertex such that $G$ becomes the disjoint union of $G'$ and $G''$. 
Recall that the graph reduced to a single vertex (without any edge) has
no isolated vertex.
\end{definition}

We directly see that $\Delta$ has for possible coefficients $0$ and
$1$, that it is coassociative and cocomutative, and that its primitive
elements are connected graphs. Here are some examples:

\begin{displaymath}
\Delta(\mathcal{M}_{{\newcommand{\nodea}{\node[draw,circle,scale=0.2,fill] (a) {$$};}
\newcommand{\nodeb}{\node[draw,circle,scale=0.2,fill] (b) {$$};}
\newcommand{\nodec}{\node[draw,circle,scale=0.2,fill] (c) {$$};}
\begin{tikzpicture}[auto]
\matrix[column sep=.0cm, row sep=.1cm,ampersand replacement=\&]{
                  \& \nodeb \&        \\ 
         \nodea   \&        \& \nodec \\ 
};
\path[thick] (a) edge (b) 
             (b) edge (c)
             (c) edge (a);
\end{tikzpicture}}}) = 
\mathcal{M}_{{\newcommand{\nodea}{\node[draw,circle,scale=0.2,fill] (a) {$$};}
\newcommand{\nodeb}{\node[draw,circle,scale=0.2,fill] (b) {$$};}
\newcommand{\nodec}{\node[draw,circle,scale=0.2,fill] (c) {$$};}
\begin{tikzpicture}[auto]
\matrix[column sep=.0cm, row sep=.1cm,ampersand replacement=\&]{
                  \& \nodeb \&        \\ 
         \nodea   \&        \& \nodec \\ 
};
\path[thick] (a) edge (b) 
             (b) edge (c)
             (c) edge (a);
\end{tikzpicture}}} \otimes 1 + 1 \otimes
\mathcal{M}_{{\newcommand{\nodea}{\node[draw,circle,scale=0.2,fill] (a) {$$};}
\newcommand{\nodeb}{\node[draw,circle,scale=0.2,fill] (b) {$$};}
\newcommand{\nodec}{\node[draw,circle,scale=0.2,fill] (c) {$$};}
\begin{tikzpicture}[auto]
\matrix[column sep=.0cm, row sep=.1cm,ampersand replacement=\&]{
                  \& \nodeb \&        \\ 
         \nodea   \&        \& \nodec \\ 
};
\path[thick] (a) edge (b) 
             (b) edge (c)
             (c) edge (a);
\end{tikzpicture}}}
\end{displaymath}

\begin{displaymath}
\Delta(\mathcal{M}_{{\newcommand{\nodea}{\node[draw,circle,scale=0.2,fill] (a) {$$};}
\newcommand{\nodeb}{\node[draw,circle,scale=0.2,fill] (b) {$$};}
\newcommand{\nodec}{\node[draw,circle,scale=0.2,fill] (c) {$$};}
\newcommand{\noded}{\node[draw,circle,scale=0.2,fill] (d) {$$};}
\begin{tikzpicture}[auto]
\matrix[column sep=.0cm, row sep=.1cm,ampersand replacement=\&]{
                  \& \nodeb \&        \& \noded \\ 
         \nodea   \&        \& \nodec \&        \\ 
};
\path[thick] (a) edge (b) 
             (c) edge (d);
\end{tikzpicture}}}) =
\mathcal{M}_{{\newcommand{\nodea}{\node[draw,circle,scale=0.2,fill] (a) {$$};}
\newcommand{\nodeb}{\node[draw,circle,scale=0.2,fill] (b) {$$};}
\newcommand{\nodec}{\node[draw,circle,scale=0.2,fill] (c) {$$};}
\newcommand{\noded}{\node[draw,circle,scale=0.2,fill] (d) {$$};}
\begin{tikzpicture}[auto]
\matrix[column sep=.0cm, row sep=.1cm,ampersand replacement=\&]{
                  \& \nodeb \&        \& \noded \\ 
         \nodea   \&        \& \nodec \&        \\ 
};
\path[thick] (a) edge (b) 
             (c) edge (d);
\end{tikzpicture}}} \otimes 1 +
\mathcal{M}_{{\newcommand{\nodea}{\node[draw,circle,scale=0.2,fill] (a) {$$};}
\newcommand{\nodeb}{\node[draw,circle,scale=0.2,fill] (b) {$$};}
\begin{tikzpicture}[auto]
\matrix[column sep=.0cm, row sep=.1cm,ampersand replacement=\&]{
                  \& \nodeb \\ 
         \nodea   \&        \\ 
};
\path[thick] (a) edge (b);
\end{tikzpicture}}} \otimes
\mathcal{M}_{{\newcommand{\nodea}{\node[draw,circle,scale=0.2,fill] (a) {$$};}
\newcommand{\nodeb}{\node[draw,circle,scale=0.2,fill] (b) {$$};}
\begin{tikzpicture}[auto]
\matrix[column sep=.0cm, row sep=.1cm,ampersand replacement=\&]{
                  \& \nodeb \\ 
         \nodea   \&        \\ 
};
\path[thick] (a) edge (b);
\end{tikzpicture}}} + 1 \otimes
\mathcal{M}_{{\newcommand{\nodea}{\node[draw,circle,scale=0.2,fill] (a) {$$};}
\newcommand{\nodeb}{\node[draw,circle,scale=0.2,fill] (b) {$$};}
\newcommand{\nodec}{\node[draw,circle,scale=0.2,fill] (c) {$$};}
\newcommand{\noded}{\node[draw,circle,scale=0.2,fill] (d) {$$};}
\begin{tikzpicture}[auto]
\matrix[column sep=.0cm, row sep=.1cm,ampersand replacement=\&]{
                  \& \nodeb \&        \& \noded \\ 
         \nodea   \&        \& \nodec \&        \\ 
};
\path[thick] (a) edge (b) 
             (c) edge (d);
\end{tikzpicture}}}
\end{displaymath}

\begin{displaymath}
\Delta(\mathcal{M}_{{\newcommand{\nodea}{\node[draw,circle,scale=0.2,fill] (a) {$$};}
\newcommand{\nodeb}{\node[draw,circle,scale=0.2,fill] (b) {$$};}
\newcommand{\nodec}{\node[draw,circle,scale=0.2,fill] (c) {$$};}
\newcommand{\noded}{\node[draw,circle,scale=0.2,fill] (d) {$$};}
\newcommand{\nodee}{\node[draw,circle,scale=0.2,fill] (e) {$$};}
\begin{tikzpicture}[auto]
\matrix[column sep=.0cm, row sep=.1cm,ampersand replacement=\&]{
                  \& \nodeb \&        \& \noded \&        \\ 
         \nodea   \&        \& \nodec \&        \& \nodee \\ 
};
\path[thick] (a) edge (b) 
             (b) edge (c)
             (d) edge (e);
\end{tikzpicture}}}) =
\mathcal{M}_{{\newcommand{\nodea}{\node[draw,circle,scale=0.2,fill] (a) {$$};}
\newcommand{\nodeb}{\node[draw,circle,scale=0.2,fill] (b) {$$};}
\newcommand{\nodec}{\node[draw,circle,scale=0.2,fill] (c) {$$};}
\newcommand{\noded}{\node[draw,circle,scale=0.2,fill] (d) {$$};}
\newcommand{\nodee}{\node[draw,circle,scale=0.2,fill] (e) {$$};}
\begin{tikzpicture}[auto]
\matrix[column sep=.0cm, row sep=.1cm,ampersand replacement=\&]{
                  \& \nodeb \&        \& \noded \&        \\ 
         \nodea   \&        \& \nodec \&        \& \nodee \\ 
};
\path[thick] (a) edge (b) 
             (b) edge (c)
             (d) edge (e);
\end{tikzpicture}}} \otimes 1 + 
\mathcal{M}_{{\newcommand{\nodea}{\node[draw,circle,scale=0.2,fill] (a) {$$};}
\newcommand{\nodeb}{\node[draw,circle,scale=0.2,fill] (b) {$$};}
\newcommand{\nodec}{\node[draw,circle,scale=0.2,fill] (c) {$$};}
\begin{tikzpicture}[auto]
\matrix[column sep=.0cm, row sep=.1cm,ampersand replacement=\&]{
                  \& \nodeb \&        \\ 
         \nodea   \&        \& \nodec \\ 
};
\path[thick] (a) edge (b)
             (b) edge (c);
\end{tikzpicture}}} \otimes
\mathcal{M}_{{\newcommand{\nodea}{\node[draw,circle,scale=0.2,fill] (a) {$$};}
\newcommand{\nodeb}{\node[draw,circle,scale=0.2,fill] (b) {$$};}
\begin{tikzpicture}[auto]
\matrix[column sep=.0cm, row sep=.1cm,ampersand replacement=\&]{
                  \& \nodeb \\ 
         \nodea   \&        \\ 
};
\path[thick] (a) edge (b);
\end{tikzpicture}}} +
\mathcal{M}_{{\newcommand{\nodea}{\node[draw,circle,scale=0.2,fill] (a) {$$};}
\newcommand{\nodeb}{\node[draw,circle,scale=0.2,fill] (b) {$$};}
\begin{tikzpicture}[auto]
\matrix[column sep=.0cm, row sep=.1cm,ampersand replacement=\&]{
                  \& \nodeb \\ 
         \nodea   \&        \\ 
};
\path[thick] (a) edge (b);
\end{tikzpicture}}} \otimes + 
\mathcal{M}_{{\newcommand{\nodea}{\node[draw,circle,scale=0.2,fill] (a) {$$};}
\newcommand{\nodeb}{\node[draw,circle,scale=0.2,fill] (b) {$$};}
\newcommand{\nodec}{\node[draw,circle,scale=0.2,fill] (c) {$$};}
\begin{tikzpicture}[auto]
\matrix[column sep=.0cm, row sep=.1cm,ampersand replacement=\&]{
                  \& \nodeb \&        \\ 
         \nodea   \&        \& \nodec \\ 
};
\path[thick] (a) edge (b)
             (b) edge (c);
\end{tikzpicture}}} 
+ 1 \otimes
\mathcal{M}_{{\newcommand{\nodea}{\node[draw,circle,scale=0.2,fill] (a) {$$};}
\newcommand{\nodeb}{\node[draw,circle,scale=0.2,fill] (b) {$$};}
\newcommand{\nodec}{\node[draw,circle,scale=0.2,fill] (c) {$$};}
\newcommand{\noded}{\node[draw,circle,scale=0.2,fill] (d) {$$};}
\newcommand{\nodee}{\node[draw,circle,scale=0.2,fill] (e) {$$};}
\begin{tikzpicture}[auto]
\matrix[column sep=.0cm, row sep=.1cm,ampersand replacement=\&]{
                  \& \nodeb \&        \& \noded \&        \\ 
         \nodea   \&        \& \nodec \&        \& \nodee \\ 
};
\path[thick] (a) edge (b) 
             (b) edge (c)
             (d) edge (e);
\end{tikzpicture}}}
\end{displaymath}

\begin{proposition}
  $(UGQSym, \cdot, \Delta)$ is a Hopf algebra. Moreover, it is filtered by the number of
  edges or by the number of nodes.

  \begin{proof}
    The product and coproduct are compatible by induction on the
    number of connected components of the operands. For $G_1$ and
    $G_2$ two non isomorphic connected graphs, the product
    $\mathcal{M}_{G_1} \cdot \mathcal{M}_{G_2}$ takes the form
    \begin{equation}~\label{prod}
      \mathcal{M}_{G_1} \cdot \mathcal{M}_{G_2} = \mathcal{M}_{G_1
        \sqcup G_2} + \sum c_{G_1, G_2}^H \mathcal{M}_{H} , 
    \end{equation}
    where the sum runs over some connected graphs H (note that if
    $G_1$ and $G_2$ are isomorphic, the coefficient of
    $\mathcal{M}_{G_1 \sqcup G_1}$ is $2$ instead of $1$). From
    this, as the coproduct stays simple for connected graphs, we
    easily check that $\Delta(f \cdot g) = \Delta(f) (\cdot \otimes \cdot)
    \Delta(g)$.

    Generally, Formula~(\ref{prod}) shows that any product
    $\mathcal{M}_{G_1} \cdot \mathcal{M}_{G_2}$ contains a leading
    term whose associated graph has as number of nodes the sum of the number of
    nodes of the operands, as number of edges the sum of the number of edges of
    the operands and as number of connected components the sum of the number of
    connected components of the operands. Other terms of the product
    are associated with graphs which are non trivial merges of $G_1$
    and $G_2$ and for which at least one edge has been joined (and
    thus at least two vertices from $G_1$ and $G_2$ has been merged).
  \end{proof}
\end{proposition}

Being commutative, cocommutative and filtered, the antipode $S$ is
uniquely determined by $S(1) = 1$ together with the induction formula:
\begin{displaymath}
  S(\mathcal{M}_G) = - \sum_{(G',G'') \in \Delta(\mathcal{M}_G) \atop{
    G'' \neq G}} \mathcal{M}_{G'} \cdot S(\mathcal{M}_{G''})
\end{displaymath}
Moreover $S$ satisfies $S^2 = Id$. We have naturally
$S(\mathcal{M}_G) = - \mathcal{M}_G$ for any $G$ connected graph.
Here are some non trivial examples.

\begin{displaymath}
  S(\mathcal{M}_{{
\newcommand{\nodea}{\node[draw,circle,scale=0.2,fill] (a) {$$};}
\newcommand{\nodeb}{\node[draw,circle,scale=0.2,fill] (b) {$$};}
\newcommand{\nodec}{\node[draw,circle,scale=0.2,fill] (c) {$$};}
\newcommand{\noded}{\node[draw,circle,scale=0.2,fill] (d) {$$};}
\begin{tikzpicture}[auto]
\matrix[column sep=.0cm, row sep=.1cm,ampersand replacement=\&]{
                  \& \nodeb \&        \& \noded \\ 
         \nodea   \&        \& \nodec \&        \\ 
};
\path[thick] (a) edge (b) 
             (c) edge (d);
\end{tikzpicture}}}) =
  \mathcal{M}_{{
\newcommand{\nodea}{\node[draw,circle,scale=0.2,fill] (a) {$$};}
\newcommand{\nodeb}{\node[draw,circle,scale=0.2,fill] (b) {$$};}
\newcommand{\nodec}{\node[draw,circle,scale=0.2,fill] (c) {$$};}
\newcommand{\noded}{\node[draw,circle,scale=0.2,fill] (d) {$$};}
\begin{tikzpicture}[auto]
\matrix[column sep=.0cm, row sep=.1cm,ampersand replacement=\&]{
                  \& \nodeb \&        \& \noded \\ 
         \nodea   \&        \& \nodec \&        \\ 
};
\path[thick] (a) edge (b) 
             (c) edge (d);
\end{tikzpicture}}} + 
  2 \mathcal{M}_{{
\newcommand{\nodea}{\node[draw,circle,scale=0.2,fill] (a) {$$};}
\newcommand{\nodeb}{\node[draw,circle,scale=0.2,fill] (b) {$$};}
\newcommand{\nodec}{\node[draw,circle,scale=0.2,fill] (c) {$$};}
\begin{tikzpicture}[auto]
\matrix[column sep=.0cm, row sep=.1cm,ampersand replacement=\&]{
                  \& \nodeb \&        \\ 
         \nodea   \&        \& \nodec \\ 
};
\path[thick] (a) edge (b)
             (b) edge (c);
\end{tikzpicture}}} +
  \mathcal{M}_{{
\newcommand{\nodea}{\node[draw,circle,scale=0.2,fill] (a) {$$};}
\newcommand{\nodeb}{\node[draw,circle,scale=0.2,fill] (b) {$$};}
\begin{tikzpicture}[auto]
\matrix[column sep=.0cm, row sep=.1cm,ampersand replacement=\&]{
                  \& \nodeb \\ 
         \nodea   \&        \\ 
};
\path[thick] (a) edge (b);
\end{tikzpicture}}}
\end{displaymath}

\begin{displaymath}
  S(\mathcal{M}_{{
\newcommand{\nodea}{\node[draw,circle,scale=0.2,fill] (a) {$$};}
\newcommand{\nodeb}{\node[draw,circle,scale=0.2,fill] (b) {$$};}
\newcommand{\nodec}{\node[draw,circle,scale=0.2,fill] (c) {$$};}
\newcommand{\noded}{\node[draw,circle,scale=0.2,fill] (d) {$$};}
\newcommand{\nodee}{\node[draw,circle,scale=0.2,fill] (e) {$$};}
\begin{tikzpicture}[auto]
\matrix[column sep=.0cm, row sep=.1cm,ampersand replacement=\&]{
                  \& \nodeb \&        \& \noded \&        \\ 
         \nodea   \&        \& \nodec \&        \& \nodee \\ 
};
\path[thick] (a) edge (b) 
             (b) edge (c)
             (d) edge (e);
  \end{tikzpicture}}}) =
  \mathcal{M}_{{
\newcommand{\nodea}{\node[draw,circle,scale=0.2,fill] (a) {$$};}
\newcommand{\nodeb}{\node[draw,circle,scale=0.2,fill] (b) {$$};}
\newcommand{\nodec}{\node[draw,circle,scale=0.2,fill] (c) {$$};}
\newcommand{\noded}{\node[draw,circle,scale=0.2,fill] (d) {$$};}
\newcommand{\nodee}{\node[draw,circle,scale=0.2,fill] (e) {$$};}
\begin{tikzpicture}[auto]
\matrix[column sep=.0cm, row sep=.1cm,ampersand replacement=\&]{
                  \& \nodeb \&        \& \noded \&        \\ 
         \nodea   \&        \& \nodec \&        \& \nodee \\ 
};
\path[thick] (a) edge (b) 
             (b) edge (c)
             (d) edge (e);
\end{tikzpicture}}} +
  6 \mathcal{M}_{{\newcommand{\nodea}{\node[draw,circle,scale=0.2,fill] (a) {$$};}
\newcommand{\nodeb}{\node[draw,circle,scale=0.2,fill] (b) {$$};}
\newcommand{\nodec}{\node[draw,circle,scale=0.2,fill] (c) {$$};}
\newcommand{\noded}{\node[draw,circle,scale=0.2,fill] (d) {$$};}
\begin{tikzpicture}[auto]
\matrix[column sep=.03cm, row sep=.04cm,ampersand replacement=\&]{
                  \& \noded \& \\
                  \& \nodea \& \\ 
         \nodeb   \&        \& \nodec \\ 
};
\path[thick] (a) edge (b)
             (a) edge (c)
             (a) edge (d);
\end{tikzpicture}}} + 
  6 \mathcal{M}_{{
\newcommand{\nodea}{\node[draw,circle,scale=0.2,fill] (a) {$$};}
\newcommand{\nodeb}{\node[draw,circle,scale=0.2,fill] (b) {$$};}
\newcommand{\nodec}{\node[draw,circle,scale=0.2,fill] (c) {$$};}
\begin{tikzpicture}[auto]
\matrix[column sep=.0cm, row sep=.1cm,ampersand replacement=\&]{
                  \& \nodeb \&        \\ 
         \nodea   \&        \& \nodec \\ 
};
\path[thick] (a) edge (b)
             (b) edge (c)
             (a) edge (c);
\end{tikzpicture}}} +
  4 \mathcal{M}_{{
\newcommand{\nodea}{\node[draw,circle,scale=0.2,fill] (a) {$$};}
\newcommand{\nodeb}{\node[draw,circle,scale=0.2,fill] (b) {$$};}
\newcommand{\nodec}{\node[draw,circle,scale=0.2,fill] (c) {$$};}
\newcommand{\noded}{\node[draw,circle,scale=0.2,fill] (d) {$$};}
\begin{tikzpicture}[auto]
\matrix[column sep=.0cm, row sep=.1cm,ampersand replacement=\&]{
                  \& \nodeb \&        \& \noded \\ 
         \nodea   \&        \& \nodec \&        \\ 
};
\path[thick] (a) edge (b)
             (b) edge (c)
             (c) edge (d);
\end{tikzpicture}}} +
  4 \mathcal{M}_{{
\newcommand{\nodea}{\node[draw,circle,scale=0.2,fill] (a) {$$};}
\newcommand{\nodeb}{\node[draw,circle,scale=0.2,fill] (b) {$$};}
\newcommand{\nodec}{\node[draw,circle,scale=0.2,fill] (c) {$$};}
\begin{tikzpicture}[auto]
\matrix[column sep=.0cm, row sep=.1cm,ampersand replacement=\&]{
                  \& \nodeb \&        \\ 
         \nodea   \&        \& \nodec \\ 
};
\path[thick] (a) edge (b)
             (b) edge (c);
\end{tikzpicture}}}
\end{displaymath}

\begin{displaymath}
  S(\mathcal{M}_{{
\newcommand{\nodea}{\node[draw,circle,scale=0.2,fill] (a) {$$};}
\newcommand{\nodeb}{\node[draw,circle,scale=0.2,fill] (b) {$$};}
\newcommand{\nodec}{\node[draw,circle,scale=0.2,fill] (c) {$$};}
\newcommand{\noded}{\node[draw,circle,scale=0.2,fill] (d) {$$};}
\newcommand{\nodee}{\node[draw,circle,scale=0.2,fill] (e) {$$};}
\newcommand{\nodef}{\node[draw,circle,scale=0.2,fill] (f) {$$};}
\begin{tikzpicture}[auto]
\matrix[column sep=.0cm, row sep=.1cm,ampersand replacement=\&]{
             \& \nodeb \&        \& \noded \&        \& \nodef \\ 
    \nodea   \&        \& \nodec \&        \& \nodee \&        \\ 
};
\path[thick] (a) edge (b) 
             (c) edge (d)
             (e) edge (f);
  \end{tikzpicture}}}) =
  - \mathcal{M}_{{
\newcommand{\nodea}{\node[draw,circle,scale=0.2,fill] (a) {$$};}
\newcommand{\nodeb}{\node[draw,circle,scale=0.2,fill] (b) {$$};}
\newcommand{\nodec}{\node[draw,circle,scale=0.2,fill] (c) {$$};}
\newcommand{\noded}{\node[draw,circle,scale=0.2,fill] (d) {$$};}
\newcommand{\nodee}{\node[draw,circle,scale=0.2,fill] (e) {$$};}
\newcommand{\nodef}{\node[draw,circle,scale=0.2,fill] (f) {$$};}
\begin{tikzpicture}[auto]
\matrix[column sep=.0cm, row sep=.1cm,ampersand replacement=\&]{
             \& \nodeb \&        \& \noded \&        \& \nodef \\ 
    \nodea   \&        \& \nodec \&        \& \nodee \&        \\ 
};
\path[thick] (a) edge (b) 
             (c) edge (d)
             (e) edge (f);
  \end{tikzpicture}}} - 
  2 \mathcal{M}_{{
\newcommand{\nodea}{\node[draw,circle,scale=0.2,fill] (a) {$$};}
\newcommand{\nodeb}{\node[draw,circle,scale=0.2,fill] (b) {$$};}
\newcommand{\nodec}{\node[draw,circle,scale=0.2,fill] (c) {$$};}
\newcommand{\noded}{\node[draw,circle,scale=0.2,fill] (d) {$$};}
\newcommand{\nodee}{\node[draw,circle,scale=0.2,fill] (e) {$$};}
\begin{tikzpicture}[auto]
\matrix[column sep=.0cm, row sep=.1cm,ampersand replacement=\&]{
                  \& \nodeb \&        \& \noded \&        \\ 
         \nodea   \&        \& \nodec \&        \& \nodee \\ 
};
\path[thick] (a) edge (b) 
             (b) edge (c)
             (d) edge (e);
\end{tikzpicture}}} -
  6 \mathcal{M}_{{\newcommand{\nodea}{\node[draw,circle,scale=0.2,fill] (a) {$$};}
\newcommand{\nodeb}{\node[draw,circle,scale=0.2,fill] (b) {$$};}
\newcommand{\nodec}{\node[draw,circle,scale=0.2,fill] (c) {$$};}
\newcommand{\noded}{\node[draw,circle,scale=0.2,fill] (d) {$$};}
\begin{tikzpicture}[auto]
\matrix[column sep=.03cm, row sep=.04cm,ampersand replacement=\&]{
                  \& \noded \& \\
                  \& \nodea \& \\ 
         \nodeb   \&        \& \nodec \\ 
};
\path[thick] (a) edge (b)
             (a) edge (c)
             (a) edge (d);
\end{tikzpicture}}} - 
  6 \mathcal{M}_{{
\newcommand{\nodea}{\node[draw,circle,scale=0.2,fill] (a) {$$};}
\newcommand{\nodeb}{\node[draw,circle,scale=0.2,fill] (b) {$$};}
\newcommand{\nodec}{\node[draw,circle,scale=0.2,fill] (c) {$$};}
\begin{tikzpicture}[auto]
\matrix[column sep=.0cm, row sep=.1cm,ampersand replacement=\&]{
                  \& \nodeb \&        \\ 
         \nodea   \&        \& \nodec \\ 
};
\path[thick] (a) edge (b)
             (b) edge (c)
             (a) edge (c);
\end{tikzpicture}}} -
  4 \mathcal{M}_{{
\newcommand{\nodea}{\node[draw,circle,scale=0.2,fill] (a) {$$};}
\newcommand{\nodeb}{\node[draw,circle,scale=0.2,fill] (b) {$$};}
\newcommand{\nodec}{\node[draw,circle,scale=0.2,fill] (c) {$$};}
\newcommand{\noded}{\node[draw,circle,scale=0.2,fill] (d) {$$};}
\begin{tikzpicture}[auto]
\matrix[column sep=.0cm, row sep=.1cm,ampersand replacement=\&]{
                  \& \nodeb \&        \& \noded \\ 
         \nodea   \&        \& \nodec \&        \\ 
};
\path[thick] (a) edge (b)
             (b) edge (c)
             (c) edge (d);
\end{tikzpicture}}} -
  6 \mathcal{M}_{{
\newcommand{\nodea}{\node[draw,circle,scale=0.2,fill] (a) {$$};}
\newcommand{\nodeb}{\node[draw,circle,scale=0.2,fill] (b) {$$};}
\newcommand{\nodec}{\node[draw,circle,scale=0.2,fill] (c) {$$};}
\begin{tikzpicture}[auto]
\matrix[column sep=.0cm, row sep=.1cm,ampersand replacement=\&]{
                  \& \nodeb \&        \\ 
         \nodea   \&        \& \nodec \\ 
};
\path[thick] (a) edge (b)
             (b) edge (c);
  \end{tikzpicture}}} -
  2 \mathcal{M}_{{
\newcommand{\nodea}{\node[draw,circle,scale=0.2,fill] (a) {$$};}
\newcommand{\nodeb}{\node[draw,circle,scale=0.2,fill] (b) {$$};}
\newcommand{\nodec}{\node[draw,circle,scale=0.2,fill] (c) {$$};}
\newcommand{\noded}{\node[draw,circle,scale=0.2,fill] (d) {$$};}
\begin{tikzpicture}[auto]
\matrix[column sep=.0cm, row sep=.1cm,ampersand replacement=\&]{
                  \& \nodeb \&        \& \noded \\ 
         \nodea   \&        \& \nodec \&        \\ 
};
\path[thick] (a) edge (b)
             (c) edge (d);
\end{tikzpicture}}} -
  \mathcal{M}_{{
\newcommand{\nodea}{\node[draw,circle,scale=0.2,fill] (a) {$$};}
\newcommand{\nodeb}{\node[draw,circle,scale=0.2,fill] (b) {$$};}
\begin{tikzpicture}[auto]
\matrix[column sep=.0cm, row sep=.1cm,ampersand replacement=\&]{
                  \& \nodeb \\ 
         \nodea   \&        \\ 
};
\path[thick] (a) edge (b);
\end{tikzpicture}}}
\end{displaymath}

\section{Connection with algebraic invariant theory}~\label{quatre}

The definition of the $\mathcal{M}_G$ obviously shows that the algebra
$UGQSym$ is composed of invariant functions under the relabeling action
of the vertices. We now show that $UGQSym$ contains a sufficient
number of invariants to separate finite graphs.

The next experiment consists in keeping only variables indexed by
small indices. For a fixed integer $n$, in our triangle of variables,
we will send to $0$ all variables $x_{ij}$ for $i > n$ (we thus keep
only an upper triangular part). Applying this restriction on the
$\mathcal{M}_G$, we define polynomials $\mathcal{P}_{n,G}$ which are
sums of monomials whose support contains only $\binom{n}{2}$ different
variables and inside which variables can appear with degree $1$. These
polynomials are the central objects of~\cite{Mikko}.

In the sequel, we insist on the fact that $\mathbb{K}$
is of characteristic $0$ and that we do not keep
relations $x_{ij}^2 = x_{ij}$ on remaining variables.

Let $n$ be a positive integer. The group acting on edges of graphs
over $n$ nodes is a permutation group (isomorphic to the symmetric
group $\SG_n$), a subgroup of the symmetric group $\SG_{\binom{n}{2}}$.
The following theorem exploits the combinatorial structure of rings 
of invariants under the action of a permutation group.

\begin{theorem}
  Let $n$ be a positive integer and $G$ be the permutation group acting on
  bi-indexed variables $\x := (x_{12}, x_{13}, \dots , x_{n-1~n})$ as
  the relabeling action on graphs over $n$ nodes. There exists a finite
  family of $\frac{\binom{n}{2}!}{n!}$ polynomials $\eta_{\lambda_i}$
  which are linear combinations of higher Specht polynomials such that
  the invariant ring $\mathbb{K}[\x]^G$ under the action of $G$ can be
  decomposed as
  \begin{equation}~\label{hiro}
    \mathbb{K}[\x]^G = \bigoplus_{\lambda \vdash \binom{n}{2}}
    \bigoplus_i \eta_{\lambda_i} \mathbb{K}[e_1, e_2, \dots e_{\binom{n}{2}}],
  \end{equation}
  where the $e_k$ are the elementary symmetric polynomials in
  $\binom{n}{2}$ variables and
  each $\eta_{\lambda_i}$ is a linear combination of some higher
  Specht polynomials $F_T^S$ with $S$ ant $T$ standard Young
  tableaux of shape the partition $\lambda$. 

  \begin{proof}
    Details are available in~\cite{CAI2015}. The algorithm presented
    in this paper shows how to compute the secondary invariants
    $\eta_{\lambda_i}$ inside irreducible representations of the
    ambient symmetric group of degree $\binom{n}{2}$.
    The slicing in Equation~(\ref{hiro}) is finer than the one presented
    in classical Hironaka decomposition because the first direct sum
    runs over partitions and not over possible degrees. As we can have
    several $G$-stable subspaces at a given degree, secondary
    invariants are better partitioned here. 
  \end{proof}
\end{theorem}

\begin{remark}
  The family formed by
  $\{e_1, e_2, \dots e_{\binom{n}{2}}\} \cup \{ \eta_{\lambda_{i}} \}_{\lambda \vdash \binom{n}{2}, i}$
  can separate all pairs of orbits of $G$ (relabeling action) when acting on
  $\mathbb{C}^{\binom{n}{2}}$ (vectors of $\binom{n}{2}$ complex numbers).
\end{remark}

The family generates the whole ring of invariants under the action of
$G$ and there always exists an invariant separating two different
orbits. Here, this huge family separates non oriented multi-graphs with
a complex number labeling each edge (also symmetric complex matrices of 
size $n$ with zeros on the diagonal).

Now, we go back to simple graphs in which, for each pair of nodes,
either we have an edge (when the label of the edge is $1$) or either we
do not have the link (the label is then $0$). Our goal is now to
reduce the number of elements of the huge separating family since
simple graphs are vectors of $\binom{n}{2}$ booleans. 

\begin{proposition}~\label{2-part}
  For any partition $\lambda \vdash n$ having at least three parts and
  any pair $(S, T)$ of standard Young tableaux of shape $\lambda$,
  the evaluation of higher Specht polynomials $F_T^S$ is zero on
  vectors of the type $\{0, 1\}^n$.

  \begin{proof}
    Higher Specht polynomials are built with Young symmetrizers which
    introduce some anti-symmetries. If the tableau has
    at least three boxes in the first column, the associated higher
    Specht will be at least divisible by a Vandermonde type factor
    over at least three variables. But a Vandermonde factor in $k$
    variables needs at least $k$ different values to be non zero
    ($(x_1 - x_2)(x_1 - x_3)(x_2 - x_3)$ is non zero if $x_1, x_2$ and $x_3$
    are pairwise different).
  \end{proof}
\end{proposition}

\begin{remark}~\label{elem-eval}
  Let $v$ be a vector of booleans of length $n$ and $i(v)$ the number
  of $1$ in $v$. Then the evaluations of the elementary symmetric
  polynomials $e_k$ on $v$ depend only on $e_1(v) = i(v)$.
  Precisely, we have:
  \begin{displaymath}
    \forall~ 1 \leqslant k \leqslant n, \qquad e_k(v) = \binom{i(v)}{k} = \binom{e_1(v)}{k}.
  \end{displaymath}
\end{remark}

We deduce from Proposition~\ref{2-part} and Remark~\ref{elem-eval}
that our separating family can be largely reduced for simple
graphs. We can keep only the first elementary symmetric polynomial
(which counts the number of edges) and combinations of higher
Specht polynomials for the relabeling action associated with partitions
having two parts.

Moreover, as we just want to separate orbits (and not generate the
whole ring of invariants), we can use the fact that the action of a
permutation $\sigma$ over an higher Specht polynomial $F_T^S$ depends
only on the irreducible representation (tableau $T$) and not on the
degree scaling (tableau $S$ for the multiplicity), on such
polynomials, we have $\sigma \cdot F_T^S = F_{\sigma \cdot T}^S$~\cite{MR1437932}.  
Therefore, Specht polynomials are sufficient
and higher Spechts bring only more copies of isomorphic $G$-stable
subspaces which contribute to generate the ring of invariants but do
not contribute to separate orbits. More details are available
in~\cite{MR1437932, CAI2015}.

Finally, noticing that a tableau of size
$\binom{n}{2}$ composed of two rows has at most
$\lfloor \frac{\binom{n}{2}}{2} \rfloor$ boxes at height $2$, we have:
\begin{theorem}
  Let $H_1$ and $H_2$ two graphs over $n$ nodes. The following
  statements are equivalents:
  \begin{enumerate}[(i)]
  \item $H_1$ and $H_2$ are isomorphic.
  \item For all graphs $G$ over at most $n$ nodes and having at most
    $\lfloor\frac{\binom{n}{2}}{2} \rfloor$ edges:
    \begin{displaymath}
      \mathcal{P}_{n,G}(H_1) = \mathcal{P}_{n,G}(H_2)
    \end{displaymath}
  \item For all graphs $G$ over at most $n$ nodes and having at most
    $\lfloor\frac{\binom{n}{2}}{2} \rfloor$ edges:
    \begin{displaymath}
      \mathcal{M}_{G}(H_1) = \mathcal{M}_{G}(H_2)
    \end{displaymath}    
  \end{enumerate}

  \begin{proof}
    We take all monomials associated with tableaux with
    two rows. They correspond to graphs over $n$ nodes having as
    edges the labels in the upper row. Applying the orbit sum operator
    (Reynold's operator up to a scalar) over these monomials, we
    obtain the required polynomials $\mathcal{P}_{n,G}$, which are more
    than needed to separate all orbits of booleans vectors.
  \end{proof}
\end{theorem}

Using the multiplicative structure of the algebra $UGQSym$, we can formulate:

\begin{corollary}~\label{isomor}
  Let $H_1$ and $H_2$ be two graphs over $n$ nodes. $H_1$ and $H_2$ are isomorphic if and only if
  $\mathcal{M}_{G}(H_1) = \mathcal{M}_{G}(H_2)$ for all connected graph $G$ over $n$ nodes having at most
  $\lfloor\frac{\binom{n}{2}}{2} \rfloor$ edges.
\end{corollary}

\section{Application to graph reconstruction}~\label{cinq}

The reconstruction of finite graphs is an old problem~\cite{MR0087949, MR0120127}. 
Given a finite graph $H$ over $n \geqslant 3$ nodes, we can build
the multi-set of vertex deleted subgraphs of $H$: $\{H_1, H_2, \dots , H_n\}$ 
which is formed by all induced subgraphs of $H$ by deleting
exactly one vertex. This process forms a map. The reconstruction
conjecture investigates if this map is injective, therefore finite
graphs may be determined by their subgraphs and finite graphs would be
reconstructibles.

\begin{proposition}[Kelly's lemma]~\label{Kelly}
  Let $n$ a positive integer, $H$ a graph over $n$ nodes and $G$ a
  graph over $r < n$ nodes. Let $\{H_1, H_2, \dots , H_n\}$ the (possibly
  multi-)set of graphs obtained from $H$ when deleting a single vertex.
  Thus, the number of ways of embedding $G$ in $H$ can be deduced from 
  the number of embedding $G$ in each $H_i$:
  \begin{displaymath}
    \mathcal{P}_{n,G}(H) = (n - r) \sum_{i=1}^n \mathcal{P}_{n,G}(H_i)
  \end{displaymath}
\end{proposition}

\begin{proposition}
  Finite graphs over $n$ nodes are reconstructible if the values of
  $\mathcal{P}_{n,G}$ for $G$ finite connected graph over $n$ nodes with at most
  $\lfloor \frac{\binom{n}{2}}{2} \rfloor$ edges can be deduced from
  the values of $\mathcal{P}_{n,H}$ for $H$ finite graph over at most
  $n-1$ nodes.

  \begin{proof}
    Use Corollary~\ref{isomor} and Kelly's lemma~\ref{Kelly}.
  \end{proof}
\end{proposition}


\textbf{Case $n = 3$:}
Graphs over $3$ nodes are isomorphic if and only if evaluations coincide on the single series
$\mathcal{M}_{{\newcommand{\nodea}{\node[draw,circle,scale=0.2,fill] (a) {$$}
;}\newcommand{\nodeb}{\node[draw,circle,scale=0.2,fill] (b) {$$}
;}\begin{tikzpicture}[auto]
\matrix[column sep=.0cm, row sep=.1cm,ampersand replacement=\&]{
                  \& \nodeb  \\ 
         \nodea   \&         \\ 
};
\path[thick] (a) edge (b);
\end{tikzpicture}}}$ since $\lfloor \frac{\binom{3}{2}}{2} \rfloor = \lfloor \frac{3}{2} \rfloor = 1$. 
Notice that the following relations:


\begin{displaymath}
\binom{\mathcal{M}_{{\newcommand{\nodea}{\node[draw,circle,scale=0.2,fill] (a) {$$}
;}\newcommand{\nodeb}{\node[draw,circle,scale=0.2,fill] (b) {$$}
;}\begin{tikzpicture}[auto]
\matrix[column sep=.0cm, row sep=.1cm,ampersand replacement=\&]{
                  \& \nodeb  \\ 
         \nodea   \&         \\ 
};
\path[thick] (a) edge (b);
\end{tikzpicture}}}}{2} =
\mathcal{M}_{{\newcommand{\nodea}{\node[draw,circle,scale=0.2,fill] (a) {$$};}
\newcommand{\nodeb}{\node[draw,circle,scale=0.2,fill] (b) {$$};}
\newcommand{\nodec}{\node[draw,circle,scale=0.2,fill] (c) {$$};}
\begin{tikzpicture}[auto]
\matrix[column sep=.0cm, row sep=.1cm,ampersand replacement=\&]{
                  \& \nodeb \&  \\ 
         \nodea   \&        \& \nodec \\ 
};
\path[thick] (a) edge (b) 
             (b) edge (c);
\end{tikzpicture}}} + 
\mathcal{M}_{{\newcommand{\nodea}{\node[draw,circle,scale=0.2,fill] (a) {$$};}
\newcommand{\nodeb}{\node[draw,circle,scale=0.2,fill] (b) {$$};}
\newcommand{\nodec}{\node[draw,circle,scale=0.2,fill] (c) {$$};}
\newcommand{\noded}{\node[draw,circle,scale=0.2,fill] (d) {$$};}
\begin{tikzpicture}[auto]
\matrix[column sep=.0cm, row sep=.1cm,ampersand replacement=\&]{
                  \& \nodeb \&        \& \noded \\ 
         \nodea   \&        \& \nodec \&   \\ 
};
\path[thick] (a) edge (b) 
             (c) edge (d);
\end{tikzpicture}}}  
\end{displaymath}

\begin{equation}~\label{binom3}
\binom{\mathcal{M}_{{\newcommand{\nodea}{\node[draw,circle,scale=0.2,fill] (a) {$$}
;}\newcommand{\nodeb}{\node[draw,circle,scale=0.2,fill] (b) {$$}
;}\begin{tikzpicture}[auto]
\matrix[column sep=.0cm, row sep=.1cm,ampersand replacement=\&]{
                  \& \nodeb  \\ 
         \nodea   \&         \\ 
};
\path[thick] (a) edge (b);
\end{tikzpicture}}}}{3} =
\mathcal{M}_{{\newcommand{\nodea}{\node[draw,circle,scale=0.2,fill] (a) {$$};}
\newcommand{\nodeb}{\node[draw,circle,scale=0.2,fill] (b) {$$};}
\newcommand{\nodec}{\node[draw,circle,scale=0.2,fill] (c) {$$};}
\begin{tikzpicture}[auto]
\matrix[column sep=.0cm, row sep=.1cm,ampersand replacement=\&]{
                  \& \nodeb \&  \\ 
         \nodea   \&        \& \nodec \\ 
};
\path[thick] (a) edge (b) 
(b) edge (c)
(a) edge (c);
\end{tikzpicture}}} + 
\mathcal{M}_{{\newcommand{\nodea}{\node[draw,circle,scale=0.2,fill] (a) {$$};}
\newcommand{\nodeb}{\node[draw,circle,scale=0.2,fill] (b) {$$};}
\newcommand{\nodec}{\node[draw,circle,scale=0.2,fill] (c) {$$};}
\newcommand{\noded}{\node[draw,circle,scale=0.2,fill] (d) {$$};}
\begin{tikzpicture}[auto]
\matrix[column sep=.03cm, row sep=.04cm,ampersand replacement=\&]{
                  \& \noded \& \\
                  \& \nodea \& \\ 
         \nodeb   \&        \& \nodec \\ 
};
\path[thick] (a) edge (b)
             (a) edge (c)
             (a) edge (d);
\end{tikzpicture}}} + 
\mathcal{M}_{{\newcommand{\nodea}{\node[draw,circle,scale=0.2,fill] (a) {$$};}
\newcommand{\nodeb}{\node[draw,circle,scale=0.2,fill] (b) {$$};}
\newcommand{\nodec}{\node[draw,circle,scale=0.2,fill] (c) {$$};}
\newcommand{\noded}{\node[draw,circle,scale=0.2,fill] (d) {$$};}
\begin{tikzpicture}[auto]
\matrix[column sep=.0cm, row sep=.1cm,ampersand replacement=\&]{
                  \& \nodeb \&        \& \noded \\ 
         \nodea   \&        \& \nodec \& \\ 
};
\path[thick] (a) edge (b) 
             (b) edge (c)
             (c) edge (d);
\end{tikzpicture}}} + 
\mathcal{M}_{{\newcommand{\nodea}{\node[draw,circle,scale=0.2,fill] (a) {$$};}
\newcommand{\nodeb}{\node[draw,circle,scale=0.2,fill] (b) {$$};}
\newcommand{\nodec}{\node[draw,circle,scale=0.2,fill] (c) {$$};}
\newcommand{\noded}{\node[draw,circle,scale=0.2,fill] (d) {$$};}
\newcommand{\nodee}{\node[draw,circle,scale=0.2,fill] (e) {$$};}
\begin{tikzpicture}[auto]
\matrix[column sep=.0cm, row sep=.1cm,ampersand replacement=\&]{
                  \& \nodeb \&        \& \noded \& \\ 
         \nodea   \&        \& \nodec \&        \& \nodee \\ 
};
\path[thick] (a) edge (b) 
             (b) edge (c)
             (d) edge (e);
\end{tikzpicture}}} +
\mathcal{M}_{{\newcommand{\nodea}{\node[draw,circle,scale=0.2,fill] (a) {$$};}
\newcommand{\nodeb}{\node[draw,circle,scale=0.2,fill] (b) {$$};}
\newcommand{\nodec}{\node[draw,circle,scale=0.2,fill] (c) {$$};}
\newcommand{\noded}{\node[draw,circle,scale=0.2,fill] (d) {$$};}
\newcommand{\nodee}{\node[draw,circle,scale=0.2,fill] (e) {$$};}
\newcommand{\nodef}{\node[draw,circle,scale=0.2,fill] (f) {$$};}
\begin{tikzpicture}[auto]
\matrix[column sep=.0cm, row sep=.1cm,ampersand replacement=\&]{
                  \& \nodeb \&        \& \noded \&        \& \nodef \\ 
         \nodea   \&        \& \nodec \&        \& \nodee \&        \\ 
};
\path[thick] (a) edge (b) 
             (c) edge (d)
             (e) edge (f);
\end{tikzpicture}}}  
\end{equation}

become for the finite case for three nodes:

\begin{displaymath}
\binom{\mathcal{P}_{{\newcommand{\nodea}{\node[draw,circle,scale=0.2,fill] (a) {$$};}
\newcommand{\nodeb}{\node[draw,circle,scale=0.2,fill] (b) {$$};}
\begin{tikzpicture}[auto]
\matrix[column sep=.0cm, row sep=.1cm,ampersand replacement=\&]{
                 \& \nodeb  \\ 
         \nodea  \&         \\ 
};
\path[thick] (a) edge (b);
\end{tikzpicture}}}}{2} =
\mathcal{P}_{{\newcommand{\nodea}{\node[draw,circle,scale=0.2,fill] (a) {$$};}
\newcommand{\nodeb}{\node[draw,circle,scale=0.2,fill] (b) {$$};}
\newcommand{\nodec}{\node[draw,circle,scale=0.2,fill] (c) {$$};}
\begin{tikzpicture}[auto]
\matrix[column sep=.0cm, row sep=.1cm,ampersand replacement=\&]{
                  \& \nodeb \&  \\ 
         \nodea   \&        \& \nodec \\ 
};
\path[thick] (a) edge (b) 
             (b) edge (c);
\end{tikzpicture}}} \qquad \qquad
\binom{\mathcal{P}_{{\newcommand{\nodea}{\node[draw,circle,scale=0.2,fill] (a) {$$}
;}\newcommand{\nodeb}{\node[draw,circle,scale=0.2,fill] (b) {$$}
;}\begin{tikzpicture}[auto]
\matrix[column sep=.0cm, row sep=.1cm,ampersand replacement=\&]{
                  \& \nodeb  \\ 
         \nodea   \&         \\ 
};
\path[thick] (a) edge (b);
\end{tikzpicture}}}}{3} =
\mathcal{P}_{{\newcommand{\nodea}{\node[draw,circle,scale=0.2,fill] (a) {$$};}
\newcommand{\nodeb}{\node[draw,circle,scale=0.2,fill] (b) {$$};}
\newcommand{\nodec}{\node[draw,circle,scale=0.2,fill] (c) {$$};}
\begin{tikzpicture}[auto]
\matrix[column sep=.0cm, row sep=.1cm,ampersand replacement=\&]{
                  \& \nodeb \&  \\ 
         \nodea   \&        \& \nodec \\ 
};
\path[thick] (a) edge (b) 
(b) edge (c)
(a) edge (c);
\end{tikzpicture}}}
\end{displaymath}

Since other patterns have as support more than $3$ nodes, they vanish
in this finite case.

\textbf{Case $n = 4$:}
As $\lfloor \frac{\binom{4}{2}}{2} \rfloor = 3$, graphs
over $4$ nodes are reconstructibles if 
$\mathcal{P}_{{\newcommand{\nodea}{\node[draw,circle,scale=0.2,fill] (a) {$$};}
\newcommand{\nodeb}{\node[draw,circle,scale=0.2,fill] (b) {$$};}
\newcommand{\nodec}{\node[draw,circle,scale=0.2,fill] (c) {$$};}
\newcommand{\noded}{\node[draw,circle,scale=0.2,fill] (d) {$$};}
\begin{tikzpicture}[auto]
\matrix[column sep=.03cm, row sep=.04cm,ampersand replacement=\&]{
                  \& \noded \& \\
                  \& \nodea \& \\ 
         \nodeb   \&        \& \nodec \\ 
};
\path[thick] (a) edge (b)
             (a) edge (c)
             (a) edge (d);
\end{tikzpicture}}}$ and
$\mathcal{P}_{{\newcommand{\nodea}{\node[draw,circle,scale=0.2,fill] (a) {$$};}
\newcommand{\nodeb}{\node[draw,circle,scale=0.2,fill] (b) {$$};}
\newcommand{\nodec}{\node[draw,circle,scale=0.2,fill] (c) {$$};}
\newcommand{\noded}{\node[draw,circle,scale=0.2,fill] (d) {$$};}
\begin{tikzpicture}[auto]
\matrix[column sep=.0cm, row sep=.1cm,ampersand replacement=\&]{
                  \& \nodeb \&        \& \noded \\ 
         \nodea   \&        \& \nodec \& \\ 
};
\path[thick] (a) edge (b) 
             (b) edge (c)
             (c) edge (d);
\end{tikzpicture}}}$ 
can be deduced from counting function of graphs over at most $3$ nodes. 
These two functions correspond to the only two connected graphs over $4$ nodes
having at most $3$ edges.

Equation~(\ref{binom3}) and the product~(\ref{prod1}) give for graphs over four nodes the identities:
\begin{displaymath}
\left\{ \begin{array}{ccc}
\binom{\mathcal{P}_{{\newcommand{\nodea}{\node[draw,circle,scale=0.2,fill] (a) {$$}
;}\newcommand{\nodeb}{\node[draw,circle,scale=0.2,fill] (b) {$$}
;}\begin{tikzpicture}[auto]
\matrix[column sep=.0cm, row sep=.1cm,ampersand replacement=\&]{
                  \& \nodeb  \\ 
         \nodea   \&         \\ 
};
\path[thick] (a) edge (b);
\end{tikzpicture}}}}{3} -
\mathcal{P}_{{\newcommand{\nodea}{\node[draw,circle,scale=0.2,fill] (a) {$$};}
\newcommand{\nodeb}{\node[draw,circle,scale=0.2,fill] (b) {$$};}
\newcommand{\nodec}{\node[draw,circle,scale=0.2,fill] (c) {$$};}
\begin{tikzpicture}[auto]
\matrix[column sep=.0cm, row sep=.1cm,ampersand replacement=\&]{
                  \& \nodeb \&  \\ 
         \nodea   \&        \& \nodec \\ 
};
\path[thick] (a) edge (b) 
(b) edge (c)
(a) edge (c);
\end{tikzpicture}}} & = &
\mathcal{P}_{{\newcommand{\nodea}{\node[draw,circle,scale=0.2,fill] (a) {$$};}
\newcommand{\nodeb}{\node[draw,circle,scale=0.2,fill] (b) {$$};}
\newcommand{\nodec}{\node[draw,circle,scale=0.2,fill] (c) {$$};}
\newcommand{\noded}{\node[draw,circle,scale=0.2,fill] (d) {$$};}
\begin{tikzpicture}[auto]
\matrix[column sep=.03cm, row sep=.04cm,ampersand replacement=\&]{
                  \& \noded \& \\
                  \& \nodea \& \\ 
         \nodeb   \&        \& \nodec \\ 
};
\path[thick] (a) edge (b)
             (a) edge (c)
             (a) edge (d);
\end{tikzpicture}}} + 
\mathcal{P}_{{\newcommand{\nodea}{\node[draw,circle,scale=0.2,fill] (a) {$$};}
\newcommand{\nodeb}{\node[draw,circle,scale=0.2,fill] (b) {$$};}
\newcommand{\nodec}{\node[draw,circle,scale=0.2,fill] (c) {$$};}
\newcommand{\noded}{\node[draw,circle,scale=0.2,fill] (d) {$$};}
\begin{tikzpicture}[auto]
\matrix[column sep=.0cm, row sep=.1cm,ampersand replacement=\&]{
                  \& \nodeb \&        \& \noded \\ 
         \nodea   \&        \& \nodec \& \\ 
};
\path[thick] (a) edge (b) 
             (b) edge (c)
             (c) edge (d);
\end{tikzpicture}}}  \\
\mathcal{P}_{{\newcommand{\nodea}{\node[draw,circle,scale=0.2,fill] (a) {$$};}
\newcommand{\nodeb}{\node[draw,circle,scale=0.2,fill] (b) {$$};}
\newcommand{\nodec}{\node[draw,circle,scale=0.2,fill] (c) {$$};}
\begin{tikzpicture}[auto]
\matrix[column sep=.0cm, row sep=.1cm,ampersand replacement=\&]{
                  \& \nodeb \&  \\ 
         \nodea   \&        \& \nodec \\ 
};
\path[thick] (a) edge (b) 
             (b) edge (c);
\end{tikzpicture}}} \cdot (
\mathcal{P}_{{\newcommand{\nodea}{\node[draw,circle,scale=0.2,fill] (a) {$$};}
\newcommand{\nodeb}{\node[draw,circle,scale=0.2,fill] (b) {$$};}
\begin{tikzpicture}[auto]
\matrix[column sep=.0cm, row sep=.1cm,ampersand replacement=\&]{
                  \& \nodeb \\ 
         \nodea   \&        \\ 
};
\path[thick] (a) edge (b);
\end{tikzpicture}}} - 2 ) - 
3 ~ \mathcal{P}_{{\newcommand{\nodea}{\node[draw,circle,scale=0.2,fill] (a) {$$};}
\newcommand{\nodeb}{\node[draw,circle,scale=0.2,fill] (b) {$$};}
\newcommand{\nodec}{\node[draw,circle,scale=0.2,fill] (c) {$$};}
\begin{tikzpicture}[auto]
\matrix[column sep=.0cm, row sep=.1cm,ampersand replacement=\&]{
                  \& \nodeb \& \\ 
         \nodea   \&        \& \nodec \\ 
};
\path[thick] (a) edge (b)
             (b) edge (c)
             (c) edge (a);
\end{tikzpicture}}} & = &
3 ~ \mathcal{P}_{{\newcommand{\nodea}{\node[draw,circle,scale=0.2,fill] (a) {$$};}
\newcommand{\nodeb}{\node[draw,circle,scale=0.2,fill] (b) {$$};}
\newcommand{\nodec}{\node[draw,circle,scale=0.2,fill] (c) {$$};}
\newcommand{\noded}{\node[draw,circle,scale=0.2,fill] (d) {$$};}
\begin{tikzpicture}[auto]
\matrix[column sep=.03cm, row sep=.04cm,ampersand replacement=\&]{
                  \& \noded \& \\
                  \& \nodea \& \\ 
         \nodeb   \&        \& \nodec \\ 
};
\path[thick] (a) edge (b)
             (a) edge (c)
             (a) edge (d);
\end{tikzpicture}}} +
2 ~ \mathcal{P}_{{\newcommand{\nodea}{\node[draw,circle,scale=0.2,fill] (a) {$$};}
\newcommand{\nodeb}{\node[draw,circle,scale=0.2,fill] (b) {$$};}
\newcommand{\nodec}{\node[draw,circle,scale=0.2,fill] (c) {$$};}
\newcommand{\noded}{\node[draw,circle,scale=0.2,fill] (d) {$$};}
\begin{tikzpicture}[auto]
\matrix[column sep=.0cm, row sep=.1cm,ampersand replacement=\&]{
                  \& \nodeb \&        \& \noded \\ 
         \nodea   \&        \& \nodec \& \\ 
};
\path[thick] (a) edge (b) 
             (b) edge (c)
             (c) edge (d);
\end{tikzpicture}}}
\end{array}\right.
\end{displaymath}

\begin{conjecture}
  Let $n \geqslant 3$ and $G$ a connected graph over $n$ nodes having at most 
  $\lfloor \frac{\binom{n}{2}}{2} \rfloor$ edges, then $\mathcal{P}_{n, G}$ is a polynomial
  over some $\mathcal{P}_{n, H}$ where $H$ are graphs over at most $n-1$ nodes.
\end{conjecture}

If this conjecture is true, it would imply Ulam's conjecture. Note
also that this conjecture will not give information about the
complexity of graph isomorphism problem.

\vspace{0.05cm}

Here is the matrix $\mathcal{M}_{G}(H)$ for $(G,H)$ the $23$ smallest
graphs without isolated nodes. 
\begin{scriptsize}
\begin{displaymath}
\newcommand{\nodea}{\node[draw,circle,scale=0.2,fill] (a) {$$};}
\newcommand{\nodeb}{\node[draw,circle,scale=0.2,fill] (b) {$$};}
\newcommand{\nodec}{\node[draw,circle,scale=0.2,fill] (c) {$$};}
\newcommand{\noded}{\node[draw,circle,scale=0.2,fill] (d) {$$};}
\newcommand{\nodee}{\node[draw,circle,scale=0.2,fill] (e) {$$};}
\begin{array}{c|c|cc|ccccccc|ccccccccccccc}
\hspace{-0.25cm}\vspace{-0.1cm} & \hspace{-0.25cm} \begin{tikzpicture}[auto]
\matrix[column sep=.0cm, row sep=.1cm,ampersand replacement=\&]{
                  \& \nodeb \\ 
         \nodea   \&        \\ 
};
\path (a) edge (b);
\end{tikzpicture} \hspace{-0.1cm} & \hspace{-0.1cm}
\begin{tikzpicture}[auto]
\matrix[column sep=.0cm, row sep=.1cm,ampersand replacement=\&]{
                  \& \nodeb \&        \\ 
         \nodea   \&        \& \nodec \\ 
};
\path (a) edge (b)
             (b) edge (c);
\end{tikzpicture} \hspace{-0.1cm} & \hspace{-0.1cm}
\begin{tikzpicture}[auto]
\matrix[column sep=.0cm, row sep=.1cm,ampersand replacement=\&]{
                  \& \nodeb \&        \\ 
         \nodea   \&        \& \nodec \\ 
};
\path (a) edge (b)
             (b) edge (c)
             (a) edge (c);
\end{tikzpicture} \hspace{-0.1cm} & \hspace{-0.1cm}
\begin{tikzpicture}[auto]
\matrix[column sep=.0cm, row sep=.1cm,ampersand replacement=\&]{
                  \& \nodeb \&        \& \noded \\ 
         \nodea   \&        \& \nodec \&        \\ 
};
\path (a) edge (b)
             (c) edge (d);
\end{tikzpicture} \hspace{-0.1cm} & \hspace{-0.1cm}
\begin{tikzpicture}[auto]
\matrix[column sep=.01cm, row sep=.02cm,ampersand replacement=\&]{
                  \& \& \nodeb \& \&         \\       
                  \& \& \noded \& \&         \\ 
         \nodea   \& \&        \& \& \nodec  \\ 
};
\path (a) edge (d)
             (b) edge (d)
             (c) edge (d);
\end{tikzpicture} \hspace{-0.1cm} & \hspace{-0.1cm}
\begin{tikzpicture}[auto]
\matrix[column sep=.0cm, row sep=.1cm,ampersand replacement=\&]{
                  \& \nodeb \&        \& \noded \\ 
         \nodea   \&        \& \nodec \&        \\ 
};
\path (a) edge (b)
             (b) edge (c)
             (c) edge (d);
\end{tikzpicture} \hspace{-0.1cm} & \hspace{-0.1cm}
\begin{tikzpicture}[auto]
\matrix[column sep=.02cm, row sep=.04cm,ampersand replacement=\&]{
                  \& \& \nodeb \& \&         \\       
                  \& \& \noded \& \&         \\ 
         \nodea   \& \&        \& \& \nodec  \\ 
};
\path (a) edge (d)
             (b) edge (d)
             (c) edge (d)
             (a) edge (c);
\end{tikzpicture} \hspace{-0.1cm} & \hspace{-0.1cm}
\begin{tikzpicture}[auto]
\matrix[column sep=.08cm, row sep=.1cm,ampersand replacement=\&]{
         \nodeb \&        \& \noded \\ 
         \nodea \&        \& \nodec \\ 
};
\path (a) edge (b)
             (b) edge (d)
             (c) edge (d)
             (a) edge (c);
\end{tikzpicture} \hspace{-0.1cm} & \hspace{-0.1cm}
\begin{tikzpicture}[auto]
\matrix[column sep=.08cm, row sep=.1cm,ampersand replacement=\&]{
         \nodeb \&        \& \noded \\ 
         \nodea \&        \& \nodec \\ 
};
\path (a) edge (b)
             (b) edge (d)
             (c) edge (d)
             (a) edge (c)
             (a) edge (d);
\end{tikzpicture} \hspace{-0.1cm} & \hspace{-0.1cm}
\begin{tikzpicture}[auto]
\matrix[column sep=.08cm, row sep=.1cm,ampersand replacement=\&]{
         \nodeb \&        \& \noded \\ 
         \nodea \&        \& \nodec \\ 
};
\path (a) edge (b)
             (b) edge (d)
             (c) edge (d)
             (a) edge (c)
             (a) edge (d)
             (b) edge (c);
\end{tikzpicture} \hspace{-0.15cm} & \hspace{-0.15cm}
\begin{tikzpicture}[auto]
\matrix[column sep=.0cm, row sep=.1cm,ampersand replacement=\&]{
                  \& \nodeb \&        \& \noded \&        \\ 
         \nodea   \&        \& \nodec \&        \& \nodee \\ 
};
\path (a) edge (b)
             (b) edge (c)
             (d) edge (e);
\end{tikzpicture} \hspace{-0.15cm} & \hspace{-0.15cm}
\begin{tikzpicture}[auto]
\matrix[column sep=.0cm, row sep=.1cm,ampersand replacement=\&]{
                  \& \nodeb \&        \& \noded \&        \\ 
         \nodea   \&        \& \nodec \&        \& \nodee \\ 
};
\path (a) edge (b)
             (b) edge (c)
             (a) edge (c)
             (d) edge (e);
\end{tikzpicture} \hspace{-0.15cm} & \hspace{-0.15cm}
\begin{tikzpicture}[auto]
\matrix[column sep=.015cm, row sep=.04cm,ampersand replacement=\&]{
         \nodeb   \& \&        \& \& \nodee  \\       
                  \& \& \noded \& \&         \\ 
         \nodea   \& \&        \& \& \nodec  \\ 
};
\path (a) edge (d)
             (b) edge (d)
             (c) edge (d)
             (e) edge (d);
\end{tikzpicture} \hspace{-0.15cm} & \hspace{-0.15cm}
\begin{tikzpicture}[auto]
\matrix[column sep=.015cm, row sep=.03cm,ampersand replacement=\&]{
         \nodee   \& \&        \& \&         \\
                  \& \& \nodeb \& \&         \\       
                  \& \& \noded \& \&         \\ 
         \nodea   \& \&        \& \& \nodec  \\ 
};
\path (a) edge (d)
             (b) edge (d)
             (c) edge (d)
             (e) edge (b);
\end{tikzpicture} \hspace{-0.15cm} & \hspace{-0.15cm}
\begin{tikzpicture}[auto]
\matrix[column sep=.0cm, row sep=.1cm,ampersand replacement=\&]{
                  \& \nodeb \&        \& \noded \&        \\ 
         \nodea   \&        \& \nodec \&        \& \nodee \\ 
};
\path (a) edge (b)
             (b) edge (c)
             (c) edge (d)
             (d) edge (e);
\end{tikzpicture} \hspace{-0.15cm} & \hspace{-0.15cm}
\begin{tikzpicture}[auto]
\matrix[column sep=.015cm, row sep=.03cm,ampersand replacement=\&]{
         \nodeb   \& \&        \& \& \nodee  \\       
                  \& \& \noded \& \&         \\ 
         \nodea   \& \&        \& \& \nodec  \\ 
};
\path (a) edge (d)
             (b) edge (d)
             (c) edge (d)
             (e) edge (d)
             (a) edge (c);
\end{tikzpicture} \hspace{-0.15cm} & \hspace{-0.15cm}
\begin{tikzpicture}[auto]
\matrix[column sep=.015cm, row sep=.04cm,ampersand replacement=\&]{
         \nodeb   \& \&        \& \& \nodee  \\       
         \nodea   \& \&        \& \& \nodec  \\ 
                  \& \& \noded \& \&         \\ 
};
\path (a) edge (c)
             (a) edge (d)
             (c) edge (d)
             (a) edge (b)
             (c) edge (e);
\end{tikzpicture} \hspace{-0.15cm} & \hspace{-0.15cm}
\begin{tikzpicture}[auto]
\matrix[column sep=.015cm, row sep=.03cm,ampersand replacement=\&]{
         \nodee   \& \&        \& \&         \\
                  \& \& \nodeb \& \&         \\       
                  \& \& \noded \& \&         \\ 
         \nodea   \& \&        \& \& \nodec  \\ 
};
\path (a) edge (d)
             (b) edge (d)
             (c) edge (d)
             (e) edge (b)
             (a) edge (c);
\end{tikzpicture} \hspace{-0.15cm} & \hspace{-0.15cm}
\begin{tikzpicture}[auto]
\matrix[column sep=.015cm, row sep=.04cm,ampersand replacement=\&]{
                  \& \& \nodeb \& \&         \\       
                  \& \& \noded \& \&         \\ 
         \nodea   \& \&        \& \& \nodec  \\ 
                  \& \& \nodee \& \&         \\ 
};
\path (a) edge (d)
             (a) edge (e)
             (e) edge (c)
             (c) edge (d)
             (d) edge (b);
\end{tikzpicture} \hspace{-0.15cm} & \hspace{-0.15cm}
\begin{tikzpicture}[auto]
\matrix[column sep=.01cm, row sep=.06cm,ampersand replacement=\&]{
                  \&        \& \nodeb \&        \&        \\ 
         \nodea   \&        \&        \&        \& \nodec \\ 
                  \& \nodee \&        \& \noded \&        \\ 
};
\path (a) edge (b)
             (b) edge (c)
             (c) edge (d)
             (d) edge (e)
             (e) edge (a);
\end{tikzpicture} \hspace{-0.15cm} & \hspace{-0.15cm}
\begin{tikzpicture}[auto]
\matrix[column sep=.015cm, row sep=.04cm,ampersand replacement=\&]{
                  \& \& \nodeb \& \&         \\       
                  \& \& \noded \& \&         \\ 
         \nodea   \& \&        \& \& \nodec  \\ 
                  \& \& \nodee \& \&         \\ 
};
\path (a) edge (d)
             (a) edge (e)
             (e) edge (c)
             (c) edge (d)
             (d) edge (b)
             (e) edge (d);
\end{tikzpicture} \hspace{-0.15cm} & \hspace{-0.15cm}
\begin{tikzpicture}[auto]
\matrix[column sep=.015cm, row sep=.04cm,ampersand replacement=\&]{
                  \& \& \nodeb \& \&         \\       
                  \& \& \noded \& \&         \\ 
         \nodea   \& \&        \& \& \nodec  \\ 
                  \& \& \nodee \& \&         \\ 
};
\path (a) edge (d)
             (a) edge (e)
             (e) edge (c)
             (c) edge (d)
             (d) edge (b)
             (a) edge (c);
\end{tikzpicture} \hspace{-0.15cm} & \hspace{-0.15cm}
\begin{tikzpicture}[auto]
\matrix[column sep=.015cm, row sep=.04cm,ampersand replacement=\&]{
         \nodeb   \&        \& \nodee  \\ 
                  \& \noded \&         \\ 
         \nodea   \&        \& \nodec  \\ 
};
\path (a) edge (d)
             (a) edge (c)
             (d) edge (c)
             (b) edge (d)
             (e) edge (b)
             (e) edge (d);
\end{tikzpicture} \\ \hline 
\begin{tikzpicture}[auto]
\matrix[column sep=.0cm, row sep=.1cm,ampersand replacement=\&]{
                  \& \nodeb \\ 
         \nodea   \&        \\ 
};
\path (a) edge (b);
\end{tikzpicture} & 1 & . & . & . & . & . & . & . & . & . & . & . & . & . & . & . & . & . & . & . & . & . & . \\ \hline
\vspace{-0.1cm} \begin{tikzpicture}[auto]
\matrix[column sep=.0cm, row sep=.1cm,ampersand replacement=\&]{
                  \& \nodeb \&        \\ 
         \nodea   \&        \& \nodec \\ 
};
\path (a) edge (b)
             (b) edge (c);
\end{tikzpicture} & 2 & 1 & . & . & . & . & . & . & . & . & . & . & . & . & . & . & . & . & . & . & . & . & . \\
\begin{tikzpicture}[auto]
\matrix[column sep=.0cm, row sep=.1cm,ampersand replacement=\&]{
                  \& \nodeb \&        \\ 
         \nodea   \&        \& \nodec \\ 
};
\path (a) edge (b)
             (b) edge (c)
             (a) edge (c);
\end{tikzpicture} & 3 & 3 & 1 & . & . & . & . & . & . & . & . & . & . & . & . & . & . & . & . & . & . & . & . \\ \hline
\vspace{-0.1cm} \begin{tikzpicture}[auto]
\matrix[column sep=.0cm, row sep=.1cm,ampersand replacement=\&]{
                  \& \nodeb \&        \& \noded \\ 
         \nodea   \&        \& \nodec \&        \\ 
};
\path (a) edge (b)
             (c) edge (d);
\end{tikzpicture} & 2 & . & . & 1 & . & . & . & . & . & . & . & . & . & . & . & . & . & . & . & . & . & . & . \\
\vspace{-0.1cm} \begin{tikzpicture}[auto]
\matrix[column sep=.01cm, row sep=.02cm,ampersand replacement=\&]{
                  \& \& \nodeb \& \&         \\       
                  \& \& \noded \& \&         \\ 
         \nodea   \& \&        \& \& \nodec  \\ 
};
\path (a) edge (d)
             (b) edge (d)
             (c) edge (d);
\end{tikzpicture} & 3 & 3 & . & . & 1 & . & . & . & . & . & . & . & . & . & . & . & . & . & . & . & . & . & . \\
\vspace{-0.1cm} \begin{tikzpicture}[auto]
\matrix[column sep=.0cm, row sep=.1cm,ampersand replacement=\&]{
                  \& \nodeb \&        \& \noded \\ 
         \nodea   \&        \& \nodec \&        \\ 
};
\path (a) edge (b)
             (b) edge (c)
             (c) edge (d);
\end{tikzpicture} & 3 & 2 & . & 1 & . & 1 & . & . & . & . & . & . & . & . & . & . & . & . & . & . & . & . & . \\
\vspace{-0.1cm} \begin{tikzpicture}[auto]
\matrix[column sep=.02cm, row sep=.04cm,ampersand replacement=\&]{
                  \& \& \nodeb \& \&         \\       
                  \& \& \noded \& \&         \\ 
         \nodea   \& \&        \& \& \nodec  \\ 
};
\path (a) edge (d)
             (b) edge (d)
             (c) edge (d)
             (a) edge (c);
\end{tikzpicture} & 4 & 5 & 1 & 1 & 1 & 2 & 1 & . & . & . & . & . & . & . & . & . & . & . & . & . & . & . & . \\
\vspace{-0.1cm} \begin{tikzpicture}[auto]
\matrix[column sep=.08cm, row sep=.1cm,ampersand replacement=\&]{
         \nodeb \&        \& \noded \\ 
         \nodea \&        \& \nodec \\ 
};
\path (a) edge (b)
             (b) edge (d)
             (c) edge (d)
             (a) edge (c);
\end{tikzpicture} & 4 & 4 & . & 2 & . & 4 & . & 1 & . & . & . & . & . & . & . & . & . & . & . & . & . & . & . \\
\vspace{-0.1cm} \begin{tikzpicture}[auto]
\matrix[column sep=.08cm, row sep=.1cm,ampersand replacement=\&]{
         \nodeb \&        \& \noded \\ 
         \nodea \&        \& \nodec \\ 
};
\path (a) edge (b)
             (b) edge (d)
             (c) edge (d)
             (a) edge (c)
             (a) edge (d);
\end{tikzpicture} & 5 & 8 & 2 & 2 & 2 & 6 & 4 & 1 & 1 & . & . & . & . & . & . & . & . & . & . & . & . & . & . \\
\begin{tikzpicture}[auto]
\matrix[column sep=.08cm, row sep=.1cm,ampersand replacement=\&]{
         \nodeb \&        \& \noded \\ 
         \nodea \&        \& \nodec \\ 
};
\path (a) edge (b)
             (b) edge (d)
             (c) edge (d)
             (a) edge (c)
             (a) edge (d)
             (b) edge (c);
\end{tikzpicture} & 6 & 12 & 4 & 3 & 4 & 12 & 12 & 3 & 6 & 1 & . & . & . & . & . & . & . & . & . & . & . & . & . \\ \hline
\vspace{-0.1cm} \begin{tikzpicture}[auto]
\matrix[column sep=.0cm, row sep=.1cm,ampersand replacement=\&]{
                  \& \nodeb \&        \& \noded \&        \\ 
         \nodea   \&        \& \nodec \&        \& \nodee \\ 
};
\path (a) edge (b)
             (b) edge (c)
             (d) edge (e);
\end{tikzpicture} & 3 & 1 & . & 2 & . & . & . & . & . & . & 1 & . & . & . & . & . & . & . & . & . & . & . & . \\
\vspace{-0.15cm} \begin{tikzpicture}[auto]
\matrix[column sep=.0cm, row sep=.1cm,ampersand replacement=\&]{
                  \& \nodeb \&        \& \noded \&        \\ 
         \nodea   \&        \& \nodec \&        \& \nodee \\ 
};
\path (a) edge (b)
             (b) edge (c)
             (a) edge (c)
             (d) edge (e);
\end{tikzpicture} & 4 & 3 & 1 & 3 & . & . & . & . & . & . & 3 & 1 & . & . & . & . & . & . & . & . & . & . & . \\
\vspace{-0.2cm} \begin{tikzpicture}[auto]
\matrix[column sep=.015cm, row sep=.04cm,ampersand replacement=\&]{
         \nodeb   \& \&        \& \& \nodee  \\       
                  \& \& \noded \& \&         \\ 
         \nodea   \& \&        \& \& \nodec  \\ 
};
\path (a) edge (d)
             (b) edge (d)
             (c) edge (d)
             (e) edge (d);
\end{tikzpicture} & 4 & 6 & . & . & 4 & . & . & . & . & . & . & . & 1 & . & . & . & . & . & . & . & . & . & . \\
\vspace{-0.1cm} \begin{tikzpicture}[auto]
\matrix[column sep=.015cm, row sep=.03cm,ampersand replacement=\&]{
         \nodee   \& \&        \& \&         \\
                  \& \& \nodeb \& \&         \\       
                  \& \& \noded \& \&         \\ 
         \nodea   \& \&        \& \& \nodec  \\ 
};
\path (a) edge (d)
             (b) edge (d)
             (c) edge (d)
             (e) edge (b);
\end{tikzpicture} & 4 & 4 & . & 2 & 1 & 2 & . & . & . & . & 1 & . & . & 1 & . & . & . & . & . & . & . & . & . \\
\vspace{-0.1cm} \begin{tikzpicture}[auto]
\matrix[column sep=.0cm, row sep=.1cm,ampersand replacement=\&]{
                  \& \nodeb \&        \& \noded \&        \\ 
         \nodea   \&        \& \nodec \&        \& \nodee \\ 
};
\path (a) edge (b)
             (b) edge (c)
             (c) edge (d)
             (d) edge (e);
\end{tikzpicture} & 4 & 3 & . & 3 & . & 2 & . & . & . & . & 2 & . & . & . & 1 & . & . & . & . & . & . & . & . \\
\vspace{-0.1cm} \begin{tikzpicture}[auto]
\matrix[column sep=.015cm, row sep=.03cm,ampersand replacement=\&]{
         \nodeb   \& \&        \& \& \nodee  \\       
                  \& \& \noded \& \&         \\ 
         \nodea   \& \&        \& \& \nodec  \\ 
};
\path (a) edge (d)
             (b) edge (d)
             (c) edge (d)
             (e) edge (d)
             (a) edge (c);
\end{tikzpicture} & 5 & 8 & 1 & 2 & 4 & 4 & 2 & . & . & . & 1 & . & 1 & 2 & . & 1 & . & . & . & . & . & . & . \\
\vspace{-0.2cm} \begin{tikzpicture}[auto]
\matrix[column sep=.015cm, row sep=.04cm,ampersand replacement=\&]{
         \nodeb   \& \&        \& \& \nodee  \\       
         \nodea   \& \&        \& \& \nodec  \\ 
                  \& \& \noded \& \&         \\ 
};
\path (a) edge (c)
             (a) edge (d)
             (c) edge (d)
             (a) edge (b)
             (c) edge (e);
\end{tikzpicture} & 5 & 7 & 1 & 3 & 2 & 5 & 2 & . & . & . & 2 & . & . & 2 & 1 & . & 1 & . & . & . & . & . & . \\
\vspace{-0.15cm} \begin{tikzpicture}[auto]
\matrix[column sep=.015cm, row sep=.03cm,ampersand replacement=\&]{
         \nodee   \& \&        \& \&         \\
                  \& \& \nodeb \& \&         \\       
                  \& \& \noded \& \&         \\ 
         \nodea   \& \&        \& \& \nodec  \\ 
};
\path (a) edge (d)
             (b) edge (d)
             (c) edge (d)
             (e) edge (b)
             (a) edge (c);
\end{tikzpicture} & 5 & 6 & 1 & 4 & 1 & 4 & 1 & . & . & . & 4 & 1 & . & 1 & 2 & . & . & 1 & . & . & . & . & . \\
\vspace{-0.15cm} \begin{tikzpicture}[auto]
\matrix[column sep=.015cm, row sep=.04cm,ampersand replacement=\&]{
                  \& \& \nodeb \& \&         \\       
                  \& \& \noded \& \&         \\ 
         \nodea   \& \&        \& \& \nodec  \\ 
                  \& \& \nodee \& \&         \\ 
};
\path (a) edge (d)
             (a) edge (e)
             (e) edge (c)
             (c) edge (d)
             (d) edge (b);
\end{tikzpicture} & 5 & 6 & . & 4 & 1 & 6 & . & 1 & . & . & 3 & . & . & 2 & 2 & . & . & . & 1 & . & . & . & . \\
\vspace{-0.2cm} \begin{tikzpicture}[auto]
\matrix[column sep=.01cm, row sep=.06cm,ampersand replacement=\&]{
                  \&        \& \nodeb \&        \&        \\ 
         \nodea   \&        \&        \&        \& \nodec \\ 
                  \& \nodee \&        \& \noded \&        \\ 
};
\path (a) edge (b)
             (b) edge (c)
             (c) edge (d)
             (d) edge (e)
             (e) edge (a);
\end{tikzpicture} & 5 & 5 & . & 5 & . & 5 & . & . & . & . & 5 & . & . & . & 5 & . & . & . & . & 1 & . & . & . \\
\vspace{-0.2cm} \begin{tikzpicture}[auto]
\matrix[column sep=.015cm, row sep=.04cm,ampersand replacement=\&]{
                  \& \& \nodeb \& \&         \\       
                  \& \& \noded \& \&         \\ 
         \nodea   \& \&        \& \& \nodec  \\ 
                  \& \& \nodee \& \&         \\ 
};
\path (a) edge (d)
             (a) edge (e)
             (e) edge (c)
             (c) edge (d)
             (d) edge (b)
             (e) edge (d);
\end{tikzpicture} & 6 & 11 & 2 & 4 & 5 & 10 & 6 & 1 & 1 & . & 3 & . & 1 & 5 & 2 & 2 & 2 & . & 1 & . & 1 & . & . \\
\vspace{-0.15cm} \begin{tikzpicture}[auto]
\matrix[column sep=.015cm, row sep=.04cm,ampersand replacement=\&]{
                  \& \& \nodeb \& \&         \\       
                  \& \& \noded \& \&         \\ 
         \nodea   \& \&        \& \& \nodec  \\ 
                  \& \& \nodee \& \&         \\ 
};
\path (a) edge (d)
             (a) edge (e)
             (e) edge (c)
             (c) edge (d)
             (d) edge (b)
             (a) edge (c);
\end{tikzpicture} & \mathbf{6} & \mathbf{10} & \mathbf{2} & \mathbf{5} & 3 & 10 & 5 & 1 & 1 & . & 5 & 1 & . & 4 & 4 & . & 2 & 2 & 1 & . & . & 1 & . \\
\begin{tikzpicture}[auto]
\matrix[column sep=.015cm, row sep=.04cm,ampersand replacement=\&]{
         \nodeb   \&        \& \nodee  \\ 
                  \& \noded \&         \\ 
         \nodea   \&        \& \nodec  \\ 
};
\path (a) edge (d)
             (a) edge (c)
             (d) edge (c)
             (b) edge (d)
             (e) edge (b)
             (e) edge (d);
\end{tikzpicture} & \mathbf{6} & \mathbf{10} & \mathbf{2} & \mathbf{5} & 4 & 8 & 4 & . & . & . & 6 & 2 & 1 & 4 & 4 & 2 & . & 4 & . & . & . & . & 1
\end{array}
\end{displaymath}
\end{scriptsize}
Even if completed with the remaining graphs over five nodes, the last two graphs are
the algebraically closest pair of graphs over at most five nodes.





\section*{acknowledgements}

This research was driven by computer exploration using the open-source
mathematical software \texttt{Sage}~\cite{sage} and its algebraic
combinatorics features developed by the \texttt{Sage-Combinat}
community~\cite{Sage-Combinat}, and in particular the support for Hopf
algebras developed by Nicolas M. Thiéry et al. and the enumeration
tools modulo the action of a permutation group developed by author.

\bibliographystyle{plain}
\bibliography{main}

\end{document}